\documentclass[12pt]{amsart}
\usepackage{amsmath}
\usepackage{amssymb}
\usepackage{amsthm}
\usepackage{amscd}
\usepackage{mathrsfs}

\usepackage{color}
\usepackage{hyperref}

\usepackage[all]{xy}
\usepackage[latin1]{inputenc}
\usepackage{mathpazo}
\usepackage[scaled=.95]{helvet}
\usepackage{courier}

\swapnumbers
\headsep=1cm \numberwithin{equation}{section}
\usepackage[colorinlistoftodos]{todonotes}

\newtheorem{theorem}{Theorem}[section]
\newtheorem{lemma}[theorem]{Lemma}
\newtheorem{proposition}[theorem]{Proposition}
\newtheorem{corollary}[theorem]{Corollary}

\theoremstyle{definition}
\newtheorem{definition}[theorem]{Definition}
\newtheorem{remark}[theorem]{Remark}
\newtheorem{remark and definition}[theorem]{Remark and Definition}
\newtheorem{remark and notation}[theorem]{Remark and Notation}

\usepackage{ulem,xpatch}

\newcommand\Sym{\operatorname{Sym}}

\newcommand\Ker{\operatorname{\Ker}}

\newcommand\Supp{\operatorname{Supp}}

\newcommand\Ann{\operatorname{Ann}}

\begin{document}

\title{Mixed multiplicity and Converse of Rees' theorem for modules}

\date{}

\author{M. D. Ferrari}
\address{UEM-PR, Brazil}
\email{mdsilva@uem.br}

\author{V. H. Jorge-Perez}
\address{Universidade de S{\~a}o Paulo -
ICMC, Caixa Postal 668, 13560-970, S{\~a}o Carlos-SP, Brazil}
\email{vhjperez@icmc.usp.br}

\author{L. C. Merighe}
\address{Universidade Estadual de Mato Grosso do Sul, Cidade Universit\'aria de Dourados, Caixa postal 351, 79804-970, Dourados-MS, Brazil}
\email{liliam.merighe@uems.com}

\date{\today}
\thanks{The second  author was supported by grant 2019/21181-0, S\~ao Paulo Research Foundation (FAPESP)}


\begin{abstract}
In this paper, we prove the converse of Rees'  mixed multiplicity theorem for modules, which extends the converse of the classical Rees' mixed multiplicity theorem for ideals given by Swanson - Theorem \ref{SwansonTheorem}. Specifically, we demonstrate the following result:

Let $(R,\mathfrak{m})$ be a $d$-dimensional formally equidimensional Noetherian local ring and $E_1,\dots,E_k$ be finitely generated $R$-submodules of a free $R$-module $F$ of positive rank $p$, with $x_i\in E_i$ for $i=1,\dots,k$. Consider \(S\), the symmetric algebra of \(F\), and \(I_{E_i}\), the ideal generated by the homogeneous component of degree 1 in the Rees algebra \([\mathscr{R}(E_i)]_1\). Assuming that $(x_1,\ldots,x_k)S$ and $I_{E_i}$ have the same height $k$ and the same radical, if the Buchsbaum-Rim multiplicity of $(x_1,\dots,x_k)$ and the mixed Buchsbaum-Rim multiplicity of the family $E_1,\dots,E_k$ are equal, i.e., ${\rm e_{BR}}((x_1,\dots,x_k)_{\mathfrak{p}};R_{\mathfrak{p}}) = {\rm e_{BR}}({E_1}_{\mathfrak{p}},\dots, {E_k}_{\mathfrak{p}},R_{\mathfrak{p}})$ for all prime ideals $\mathfrak{p}$ minimal over $((x_1,\ldots,x_k):_RF)$, then $(x_1,\ldots,x_k)$ is a joint reduction of $(E_1,\dots,E_k)$. 

In addition to proving this theorem, we establish several properties that relate joint reduction and mixed Buchsbaum-Rim multiplicities.

\end{abstract}
\maketitle

\section{Introduction} 
One of the main motivations for studying mixed multiplicities of ideals and modules arises from the geometric interpretation discovered in 1973 by Teissier in his paper at Cargese \cite{Teissier}. Let $R = \mathbb{C}\{z_0, z_1, \dots, z_d\}$ be the local ring of convergent power series with maximal ideal $\mathfrak{m}$. He established a significant result: if $f \in R$ represents the equation of an isolated hypersurface singularity $(X,0)$, then the Jacobian of $f$, denoted by $J(f)$, is an $\mathfrak{m}$-primary ideal. Consequently, the function $\ell_R\left(R/J(f)^u\mathfrak{m}^v\right)$ is represented by a polynomial $P(u,v)$ of total degree $d + 1$ for large values of $u$ and $v$. Moreover, the terms of total degree $d + 1$ in $P(u,v)$ can be expressed as:
$$\sum_{i=0}^{d+1}\frac{1}{(d+1-i)!\cdot i!}\mu^{(d+1-i)} u^{d +1-i}\cdot v^i,$$
where $\mu^{(i)}$ denotes the Milnor number of $X \cap E$ at $0\in \mathbb{C}^{d+1}$, and $E$ represents an $i$-dimensional affine subspace of $\mathbb{C}^{d+1}$ passing through $0\in \mathbb{C}^{d+1}$ for all $i=0,\dots, d+1$. Furthermore, within the same context, Teissier \cite{Teissier} obtained two interesting results: the first one states that the mixed multiplicities of $\mathfrak{m}$ and $J(f)$, denoted by ${\rm e}\left(\mathfrak{m}^{[d+1-i]},J(f)^{[i]}\right)$, are equal to $\mu^{(d+1-i)}$ for all $i =0,\dots,d+1$. The second result states that ${\rm e}\left(\mathfrak{m}^{[d+1-i]},J(f)^{[i]}\right)$ is equal to the Hilbert-Samuel multiplicity of $(a_1,\dots,a_{d+1-i},b_1,\dots,b_i)$, denoted by ${\rm e}\left(a_1,\dots,a_{d+1-i},b_1,\dots,b_i\right)$, where $a_1,\dots,a_{d+1-i}$ represent the defining equation of $E$ and $b_1,\dots,b_i$ are general elements in $J(f)$, for general theory, see the survey \cite[p. 535-537]{TrungVerma}). 

Since the beginning, both the Milnor number and the Hilbert-Samuel multiplicity have played crucial roles in singularity theory, commutative algebra, and algebraic geometry. Therefore, the results established by Teissier have served as a significant motivation for the development of the theory of mixed multiplicities of ideals and modules. To determine further notations and clarify some results, let us review some definitions and some fundamental properties of multiplicities.

Let $(R,\mathfrak{m})$ be a $d$-dimensional Noetherian local ring and $I_1,\ldots,I_k$ $\mathfrak{m}$-primary ideals in $R$. There exists a polynomial $P(n_1,\ldots,n_k)$ in $k$ variables with rational coefficients and total degree $d$ such that for sufficiently large $n_1,\ldots ,n_k$:
$$P\left(n_1,\ldots,n_k\right)=\ell_R\left(\frac{R}{I^{n_1}_1\cdots I^{n_k}_k}\right).$$
This polynomial $P(n_1,\ldots,n_k)$ is known as the \textit{multi-graded Hilbert polynomial} of $I_1,\ldots,I_k$. By decomposing this polynomial $P(n_1,\ldots,n_k)$ into its homogeneous parts, we can express it as:
$$\sum_{d_1+\cdots+d_k=d}\frac{1}{d_1!\cdots d_k!}
{\rm e}\left(I^{[d_1]}_1,\dots,I^{[d_k]}_k\right) n^{d_1}_1\cdots n^{d_k}_k,$$
where ${\rm e}\left(I^{[d_1]}_1,\ldots,I^{[d_k]}_k\right)\in \mathbb{Q}$ and $I^{[d_i]}_i$ indicates that each $I_i$ is listed $d_i$ times for all $i=1,\dots,k$. This number ${\rm e}\left(I^{[d_1]}_1,\dots,I^{[d_k]}_k\right)$ is called the \textit{mixed multiplicity of type $(d_1,\dots,d_k)$ with respect to $I_1,\dots,I_k$}. For more details, see \cite{SH, Rober}.

A $k$-tuple of elements of $R$ $(x_1,\ldots,x_k)$ is a \textit{joint reduction} of the $k$-tuple $(I_1,\ldots,I_k)$ if the ideal $\sum^k_{i=1}x_iI_1\cdots I_{i-1}I_{i+1} \cdots I_k$ is a reduction of $I_1\cdots I_k$, where $x_i\in I_i$ for each $i=1,\ldots,k$. It is worth noting that when $k=1$, the mixed multiplicity coincides with the Hilbert-Samuel multiplicity, denoted by ${\rm e}(I_1)$, and the definition of joint reduction aligns with the definition of reduction of an ideal. In this context, Rees in \cite{Rees1} established a fundamental result known as Rees' theorem: 

\begin{theorem}[Rees' theorem] Let $(R,\mathfrak{m})$ be a formally equidimensional Noetherian local ring and two $\mathfrak{m}$-primary ideals $(x_1,\dots,x_d)\subset I$. Then $(x_1,\dots,x_d)$ is a reduction of $I$ if and only if ${\rm e}(I) = {\rm e}(x_1,\dots,x_d)$.
\end{theorem}

Bo\"ger \cite{boger} extended this result to non-$\mathfrak{m}$-primary ideals. In addition to the results presented by Teissier, Rees and B\"oger, it also has as motivating factors for obtaining the fundamental results in the theory of mixed multiplicities. The first significant result was given by Rees \cite{Rees}, commonly known as "Rees' mixed multiplicity theorem for ideals", which establishes a connection between joint reductions and mixed multiplicities (see also \cite[Theorem 17.4.9]{SH}).

\begin{theorem}[Rees' mixed multiplicity theorem for ideals] Let $(R,\mathfrak{m})$ be a $d$-dimensional Noetherian local ring, $I_1,\ldots,I_k$ be $\mathfrak{m}$-primary ideals in $R$, $x_i\in I_i$ for all $i=1,\ldots,k$. Consider integers $d_1,\dots,d_k$ with $d_1+\cdots+d_k=d$. If $(x_1,\ldots,x_d)$ is a joint reduction of $I_1,\ldots,I_k$ , then ${\rm e}\left(x_1,\ldots,x_d\right)={\rm e}\left(I_1^{[d_1]},\ldots,I_k^{ [d_k]}\right).$
\end{theorem}

This theorem states that the ideal generated by a joint reduction has the same multiplicity as the mixed multiplicity of ideals. The second fundamental result was given by Swanson, where she showed the converse of "Rees' mixed multiplicity theorem for ideals" \cite[Theorem 17.6.1]{Swanson, SH} and further generalizes the result given by B\"oger for mixed multiplicities of ideals.

\begin{theorem}[Converse of Rees' mixed multiplicity theorem for ideals] \label{SwansonTheorem}
Let $(R,\mathfrak{m})$ be a formally equidimensional $d$-dimensional Noetherian local ring, $I_1 ,\ldots,I_k$ be ideals in $R$, $x_i\in I_i$ for $i=1,\ldots,k$. Consider integers $d_1,\dots,d_k$ with $d_1+\cdots +d_k=d$. If
${\rm e}\left((x_1,\ldots,x_k)R_{\mathfrak{p}};R_{\mathfrak{p}}\right) = {\rm e}\left(I_1R_{\mathfrak{p}} , \ldots , I_kR_{\mathfrak{p}} ;R_{\mathfrak{p}}\right)$
for all prime ideals ${\mathfrak{p}}$ minimal over $(x_1,\ldots,x_k)$, then $(x_1,\ldots,x_k)$ is a joint reduction of $(I_1,\ldots, I_k)$.
\end{theorem}

Extending the concept of multiplicity for modules over a local Noetherian ring, known as the Buchsbaum-Rim multiplicity is also well studied by many authors, for a formal definition see Section \ref{section4.2}. This multiplicity extends the notion of multiplicity for ideals, and it has been the subject of several generalizations and results. Notably, Rees' theorem and related results have been established by Buchsbaum, Rim, Kirby, Rees \cite{Kirby-Rees1}, and Katz \cite{Katz}. Moreover, Kleiman and Thorup provided an algebro-geometric interpretation of these same results in \cite{Kleiman-Thorup} and \cite{Kleiman-Thorup2}.

Katz made significant use of the Buchsbaum-Rim multiplicity to prove a reduction criterion for modules, which represents a generalization of Boger's theorem. 

Kirby-Rees, Kleiman and Thorup in \cite{Kirby-Rees1, Kleiman-Thorup, Kleiman-Thorup2} also introduced the notion of mixed multiplicities for a family $E_1, \dots, E_k$ of $R$-submodules of $R^p$ with finite colength. This multiplicity is denoted by ${\rm e_{BR}}\left(E_1^{[d_1]},\ldots,E_k^{[d_k]}\right)$, and its formal definition can be found in Section \ref{section4.3}. In this context, Bedregal-Perez \cite{Bedregal-Perez2} provided a generalization of "Rees' mixed multiplicity theorem for ideals".

\begin{theorem}[Rees' mixed multiplicity theorem for modules] \label{Perez-Bedregal}
Let $(R,\mathfrak{m})$ be a $d$-dimensional Noetherian local ring, $E_1,\dots, E_{k}$ be $R$-submodules of $F$ such that $F/E_i$ has finite length, for all $i=1, \dots, k$. Let $x_1,\dots ,x_{d+p-1}$ be any joint reduction of $E_1,\dots, E_{k}$, with each $E_i$ listed $d_i$ times. Then, for all integers $d_1,\ldots,d_k\in\mathbb{N}$ with $d_1+d_2+\cdots+d_k=d+p-1$, $
{\rm e_{BR}}\left(E_1^{[d_1]},\ldots,E_k^{[d_k]}\right) =
{\rm e_{BR}}\left(x_1,\dots ,x_{d+p-1}\right),$
where ${\rm e_{BR}}(x_1,\dots ,x_{d+p-1})$ denotes the Buchsbaum-Rim multiplicity of the $R$-submodule of $F$ generated by $x_1 ,\dots ,x_{d+p-1}$.
\end{theorem}

One of the main results of this work is to provide a generalization of B\"oger's version of the result for modules, as given by Katz, and to demonstrate the converse of "Rees' mixed multiplicity theorem for modules", the Theorem \ref{Perez-Bedregal} above.

The organization of this paper is given as follows:

In Section \ref{section2}, we will first establish some notations, definitions and basic results that will be useful for the rest of the paper.

In Section \ref{section3}, we will review some definitions and results about Buchsbaum-Rim multiplicity, mixed Buchsbaum-Rim multiplicity, and g-multiplicity. It should be noted that these invariants will be fundamental in the proofs of almost all the results in the paper.

In Section \ref{section4}, we state and prove one of the main results of the paper, called the "Risler-Teissier Theorem", for finitely generated $R$-submodules $E_i$ of $F^{e_i}$, where $e_i \in \mathbb{N}$. This result generalizes the Risler-Teissier theorem for ideals given in \cite[Proposition 2.1]{Teissier}, and also generalizes the Risler-Teissier theorem for modules given in \cite[Theorem 5.2]{Bedregal-Perez2} when $e_i=1$, for all $i$.

In Section \ref{section5}, we present fundamental results that are crucial to the proof of the main theorem of this paper.

In the last section, we establish the main result of this paper, which is the converse of Theorem \ref{Perez-Bedregal} and also a generalization of Swanson's theorem \ref{SwansonTheorem}.

\section{Setup and Background}\label{section2}

In this paper, $(R, \mathfrak{m})$ denotes a Noetherian local ring with a maximal ideal $\mathfrak{m}$, and $F$ represents a free $R$-module of positive rank $p$. Furthermore, we assume that all $R$-modules discussed in this paper are finitely generated. Let $E$ be a finitely generated submodule of $F$.  The embedding $E\subseteq F$ induces a graded $R$-algebra homomorphism between the symmetric algebras:
$${\rm Sym}(i):{\rm Sym}_R(E)\longrightarrow {\rm Sym}_R(F),$$
which maps the symmetric algebra of $E$ to the symmetric algebra of $F$. Since $F$ is a free $R$-module, the symmetric algebra $S:={\rm Sym}_R(F)$ of $F$ coincides with the polynomial ring $R[t_1,t_2,\dots,t_p]$ over $R$ with $rank(F)=p$ variables. In \cite{Simis-Ulrich-Vasconcelos}, the authors Simis-Ulrich-Vasconcelos defined the Rees algebra $\mathscr{R}(E)$ of the module $E$ as the image of this homomorphism.

The Rees algebra $\mathscr{R}(E)$ is given by:
$$\mathscr{R}(E)={\rm Im}\left({\rm Sym}(i)\right) = \bigoplus_{n\geq0} E^n\subseteq S = R[t_1,t_2,\ldots,t_p],$$
where $E^n:=[\mathscr{R}(E)]_n$ represents the homogeneous component of $\mathscr{R}(E)$ of degree $n$. In other words, $E^n$ is the $n$-th power of the image of $E$ in $\mathscr{R}(E)$. If we consider an element $h = (h_1,\dots ,h_p) \in F$, we define the element $w(h) = h_1t_1 + \cdots + h_pt_p \in S_1$. We denote by $[\mathscr{R}(E)]_1:= \{w(h) : h \in E\}=E$. In particular, $E=[\mathscr{R}(E)]_1$ is an $R$-submodule of $\mathscr{R}(E)$. The ideal of $S$ generated by $[\mathscr{R}(E)]_1$ is denoted by $I_E$.

To simplify the notation, the element $w(h)$ will be denoted simply by $h$, throughout all the paper.

Furthermore, now we establish a consistent set of notations to enhance readability and ensure clarity in our presentation. The following list outlines the notations that will be utilized throughout this paper:

\noindent \textbf{Notations:}\label{notation1} Let $R$ be a $d$-dimensional Noetherian ring, $Y$ be a variable over $R$ and $\mathfrak{p}$ a prime ideal of $R$. Let $N$ be a finitely generated $R$-module of dimension $d$:
\begin{itemize}
  \item $F$: A free $R$-module of positive rank $p$.
  \item $S = \Sym_R(F) = R[t_1, \dots, t_p]$: The symmetric algebra of $F$ over $R$.
  \item $S_{\mathfrak{p}} = \Sym_{R_{\mathfrak{p}}}F_{\mathfrak{p}}$: The localization of $S$ in $\mathfrak{p}$.
  \item $M=S\otimes_RN=\oplus_{q\geq0}M_q=\oplus_{q\geq0}F^q\otimes_RN$. Note that $M$ is a finitely generated graded $S$-module of dimension $d+p$. 
  
  \item $F[Y] = R[Y]^p$: The $R[Y]$-free module of rank $p$ generated in variable $Y$.
  \item $S[Y] = \Sym_{R[Y]}F[Y]$: The symmetric algebra of $F[Y]$ over $R[Y]$.
  \item $M[Y] = S[Y] \otimes_{R} N = \bigoplus_{q \geq 0} F[Y]^q \otimes_R N$.
  \item $E(B)$: The image of $E \subseteq R^p$ in the $B$-free module $B^p$ by the ring homomorphism $\varphi: R \to B$. For instance, if $B = R/\mathfrak{p}$, then $E(R/\mathfrak{p})$ corresponds to the image of $E$ in the $R/\mathfrak{p}$-free module $F/\mathfrak{p}F$, i.e., $E(F/\mathfrak{p}F) = (E+\mathfrak{p}F)/\mathfrak{p}F$; and  $E': =EF' =E(F/(x)) = (E+(x))/(x)$.
  \item Let ${\bf n} = (n_1, \ldots, n_k)$ and ${\bf e} = (e_1, \ldots, e_k)$ denote multi-indices. The norm $|{\bf n}|$ is represented by $|{\bf n}| = n_1 + \ldots + n_k$, the inner product of ${\bf e}$ and ${\bf n}$ is given as ${\bf e} \cdot {\bf n} = e_1n_1 + \ldots + e_kn_k$, and ${\bf n}!$ denotes the product of factorials ${\bf n}! = n_1! \cdots n_k!$. 

If ${\bf d} = (d_1, \ldots, d_k)$ is another multi-index, we denote ${\bf n}^{\bf d} = n_1^{d_1} \cdots n_k^{d_k}$.

Now, let ${\bf E} = (E_1, \ldots, E_k)$ represent the multi-index of submodules. We denote ${\bf E}^{\bf n} = E_1^{n_1} \cdots E_k^{n_k}$.

Additionally, we use the notation $\delta(i) = (\delta(i,1), \ldots, \delta(i,k))$, where $\delta(i,j) = 1$ if $i = j$ and $\delta(i,j) = 0$ otherwise.
\end{itemize}

Now, we are ready to give some definitions and results that will be used throughout the paper.

\begin{definition}[Reduction of a module] An $R$-submodule $U$ of $E$ is called a  \textit{reduction of $E$} if there exists an integer $n_0 > 0$ such that $E^{n+1} = UE^n$ for all $n \geq n_0$. A reduction $U \subseteq E$ is said to be a \textit{minimal reduction of $E$} if $U$ properly contains no further reductions of $E$.
\end{definition}

If $U$ is a reduction of $E$, then there exists at least one $R$-submodule $L$ in $U$ such that $L$ is a minimal reduction of $U$. The proof of this statement follows the same deduction as the proof of Huneke-Swanson in \cite[Theorem 8.3.5]{SH} or in Ooishi in \cite[Theorem 3.3]{ooishi}.

\begin{definition}[Integral closure of a module] Given an integer $n \geq 0$, the \textit{integral closure of the $n$-th power module $E^n$} is given as
$$\overline{E^n} = \left(\overline{\mathscr{R}(E)}^S\right)_n \subseteq S_n = F^n,$$
which represents the $n$-th homogeneous component of the integral closure $\overline{\mathscr{R}(E)}^S$ of $\mathscr{R}(E)$ in $S$.
\end{definition}
In other words, $\overline{E^n}$ is the integral closure of the ideal $(ES)^n$ of degree $n$. In particular, $\overline{E} = (\overline{ES})_1 \subseteq F$. Therefore, $\overline{E}$ consists of the elements $x \in F$ that satisfy the integral equation $x^n + c_1x^{n-1} + \cdots + c_{n-1}x + c_n = 0$ in $S$, where $n > 0$ and $c_i \in E^i$ for every $1 \leq i \leq n$. 

Next, we have an elementary fact relating reduction and integral closure.

\begin{remark} \cite[Remark 2.2]{PF}. \label{remarkideal} Suppose that $L \subseteq E$ are two $R$-submodules, and let $I_L$ and $I_E$ be the ideals of $S$. The following conditions are equivalent:
\begin{itemize}
\item[$(i)$] $L$ is a reduction of $E$.
\item[$(ii)$] $I_L$ is a reduction of $I_E$.
\item[$(iii)$] $I_E \subseteq \overline{I_L}$.
\item[$(iv)$] $E \subseteq \overline{L}$.
\end{itemize}
\end{remark}

\begin{lemma} \label{Proposition 1.1.5} Let $R$ be a ring, not necessarily Noetherian. An element $r \in F$ is in the integral closure of $E$ if and only if for every minimal prime ideal $\mathfrak{p}$ in $R$, the image of $r$ in $F/\mathfrak{p}F$ is in the integral closure of $E\left(F/\mathfrak{p}F\right)$, denoted as $(E + \mathfrak{p}F)/\mathfrak{p}F$.
 \end{lemma}
\begin{proof} By \cite[Proposition 1.1.5]{SH}, it is stated that an element $r \in S$ is in the integral closure of $I_E$ if and only if for every minimal prime ideal $\mathfrak{p}$ in $S$, the image of $r$ in $S/\mathfrak{p}$ is in the integral closure of $(I_E + \mathfrak{p})/\mathfrak{p}$. Using this result, along with Remarks \ref{remarkideal} and the persistence of the integral closure, we can conclude the desired result.
\end{proof}

The definitions of joint reduction and superficial sequence, in this context, can be given similarly to the definitions given in \cite[Definition 3.1 and 3.4]{Bedregal-Perez2}.

\begin{definition}[Joint reduction] Let $E_1,\ldots, E_{k}$ be $R$-submodules of $F$. A sequence of elements $x_1,\ldots, x_{k}$ with $x_i \in E_i$ is called a \textit{joint reduction} for $E_1,\ldots, E_{k}$ with respect to the $R$-module $N$ if there exists an integer $n\geq0$ such that for all $q \geq 0$, the following equality holds:
$$
\left[ 
\sum_{i=1}^k
x_iE_1 \cdots E_{i-1}E_{i+1}\cdots E_k
\right](E_1 \cdots E_k)^{n-1}M_q = (E_1\cdots E_k)^nM_q,
$$
where $M_q=F^q\otimes_RN$, for a positive integer $q$.

In the special case when $N=R$, we say that the sequence $x_1,\ldots, x_{k}$ is a joint reduction for $E_1,\ldots, E_{k}$.
\end{definition}

\begin{remark}\label{Remark 4.4} Notice that a sequence of elements $x_1,\dots,x_k$ with $x_i \in E_i$ is a joint reduction for $E_1,\dots,E_k$ with respect to $N$ if and only if the sequence $w(x_1),\dots,w(x_k)$ is a joint reduction for $I_{E_1},\dots,I_{E_k}$ with respect to $M=S\otimes_RN$, where $w(x_i)$ denotes the image of $x_i$ in $S\otimes_RN$. This follows from Rees' definition \cite{Rees} and the definition given by Huneke-Swanson \cite[Definition 17.1.1]{SH}. For more details, see those references.
\end{remark}


\begin{proposition}\label{Proposition 17.3.3}
Let $(R,\mathfrak{m})$ be a Noetherian local ring and $(x_1, \dots , x_k)$ be an $R$-submodule of $F$ such that $\ell_R(F/(x_1, \dots , x_k))<\infty$. Then, the sequence $(x_1, \dots , x_k)$ is a joint reduction for the submodule sequence $((x_1 ) + \mathfrak{m}^nF , \dots , (x_k) + \mathfrak{m}^nF)$ for all sufficiently large $n$.
\end{proposition}

\begin{proof}
Since $\ell_R(F/(x_1, \dots , x_k))<\infty$, there exists $n$ such that $\mathfrak{m}^nF \subseteq (x_1, \dots , x_k)$. Let
$$E = \sum_{i=1}^kx_i\left((x_1 ) + \mathfrak{m}^nF\right)\cdots \left((x_{i-1} ) + \mathfrak{m}^nF\right)\left((x_{i+1}) + \mathfrak{m}^nF\right)\cdots ((x_k ) + \mathfrak{m}^nF).$$
Then
\begin{align*}
((x_1 ) + \mathfrak{m}^nF)\cdots ((x_k) + \mathfrak{m}^nF) &= E + \mathfrak{m}^{n}F^k\\
&\subseteq E + (x_1,\dots,x_k)F^{k-1}\\ 
&\subseteq E\\
&\subseteq \left((x_1 ) + \mathfrak{m}^nF\right)\cdots \left((x_k) + \mathfrak{m}^nF\right).
\end{align*}
Therefore the equality holds, and it follows that the submodule sequence $((x_1 ) + \mathfrak{m}^nF , \dots , (x_k) + \mathfrak{m}^nF)$ forms a joint reduction.
\end{proof}

The definition of a superficial element for $E_1, \ldots, E_k$ with respect to $N$ can be stated as follows:

\begin{definition}[Superficial element] \label{defsuperficial}
Let $E_1,\ldots, E_{k}$ be $R$-submodules of $F$ and $N$ an $R$-module. An element $x\in E_1$ is a {\it superficial element for $E_1,\ldots,E_k$ with respect to $N$} if there exists a non-negative integer $c_1$ such that for all $n_1\geq c_1$ and all $n_2,\ldots,n_k\geq 0$, $q\geq 0$,
$$\left({\bf E}^{{\bf n}}M_q:_{M_{|{\bf n}|-1+q}} x\right)\cap E_1^{c_1}E_2^{n_2}\cdots E_k^{n_k}M_{n_1-1-c_1+q}=
E_1^{n_1-1}E_2^{n_2}\cdots E_k^{n_k}M_q.$$
where $M_q=F^q\otimes_RN$. 
\end{definition}

\begin{remark}\label{relasuperficial}
    Notice that $x_1 \in E_1$ is a superficial element for $E_1,\ldots, E_{k}$ with respect to $N$ in the above sense if and only if $w(x_1)\in I_{E_1}$ is a superficial element for $I_{E_1},\ldots, I_{E_{k}}$ with respect to $N$ (in the sense of \cite[Definition 17.2.1]{SH}).

A sequence of elements $x_1,\dots,x_k,$ with $x_i \in E_i$, is a {\it superficial sequence for $E_1,\dots,E_k$ with respect to $N$} if for all $i=1,\dots,k$, $x_i\in E_i$, is such that the sequence $w(x_1),\dots,w(x_k)$ is a superficial sequence for $I_{E_1},\dots,I_{E_k}$ with respect to $M=S\otimes_RN$. 
\end{remark}

\begin{proposition} (Existence of superficial elements) \cite[Proposition 3.3]{Bedregal-Perez2}\label{propo3.3BP} Let $(R,\mathfrak{m})$ be a Noetherian local ring with an infinite residue field and $E_1,\dots, E_k$ be $R$-submodules of $F$. Then there exists a superficial sequence for $E_1,\dots, E_{k}$ with respect to the $R$-module $N$. Specifically, there exists a non-empty Zariski-open subset $U$ of $E_1/\mathfrak{m}E_1$ such that for any $x \in E_1$ with image in $U$, $x$ is a superficial element for $E_1,\dots, E_k$ with respect to $N$. Moreover, if $E_1$ is not contained in the prime ideals $\mathfrak{p}_1, \dots, \mathfrak{p}_s$ of $R$. Then $x$ can be chosen to avoid the same prime ideals.
\end{proposition}

\section{ Buchsbaum-Rim multiplicity, mixed multiplicity and g-multiplicity}\label{section3}

 In this section, we review some definitions and results about Buchsbaum-Rim multiplicity, mixed Buchsbaum-Rim multiplicity and g-multiplicity, and also we give some basic properties of those multiplicities. They will be fundamental for the rest of the paper.
 
\subsection{Buchsbaum-Rim multiplicity} \label{secBR}
 
Suppose that $(R,\mathfrak{m})$ is a $d$-dimensional Noetherian local ring, $E$ is an $R$-submodule of $F$ such that $0<\ell_R(F/E)<\infty$ and Let $N$ be an $R$-module of dimension $d$. Buchsbaum-Rim, in \cite[Section 3]{buchsbaum}, showed that there exists an integer ${\rm e}_{\rm BR}(E;N)\in \mathbb{Z}$  such that the lenght $\ell_R\left(\frac{F^n\otimes_RN}{E^n\otimes_RN}\right)$ can be expressed as the polynomial of the form:
$$P_E^F(n;N)={\rm e}_{\rm BR}\left(E;N\right)\frac{n^{d+p-1}}{\left(d+p-1\right)!}+\mbox{lower terms}$$
for $n$ sufficiently large. The integer ${\rm e}_{\rm BR}(E;N)$ is called the {\it Buchsbaum-Rim multiplicity of $E$ with respect to $N$}. This multiplicity is also defined as
$${\rm e}_{\rm BR}(E;N)=\lim\limits_{n\to \infty}\left(d+p-1\right)!\frac{\ell_R\left(\frac{F^n\otimes_RN}{E^n\otimes_RN}\right)}{n^{d+p-1}}.$$
When $N = R$ we will denote ${\rm e}_{\rm BR}(E):={\rm e}_{\rm BR}(E;R)$. 

Let $E_1 \subseteq E_2$ be $R$-submodules contained in  $F$, such that $\ell_R(F/E_1) < \infty$ and $R$ is formally equidimensional ring. It is well known that $E_1$ is a reduction of $E_2$ if and only if the Buchsbaum-Rim multiplicity of $E_1$ is equal to the Buchsbaum-Rim multiplicity of $E_2$, i.e., ${\rm e}_{\rm BR}(E_1) = {\rm e}_{\rm BR}(E_2)$. For more details, see, for example, \cite{buchsbaum}, \cite{Katz}, \cite[p. 317, Corollary 16.5.7]{SH}, and \cite[Theorem 5.7]{Kleiman-Thorup}. This result generalizes Rees' theorem for ideals. More explicitly, we have the following: 

\begin{lemma}[Rees' theorem] \label{Coro16.5.7SH} Let $(R,\mathfrak{m})$ be a formally equidimensional $d$-dimensional Noetherian local ring and $E_1\subseteq E_2$ be $R$-submodules of $F$ with $\ell_R(F/E_i)<\infty$, for $i=1,2$. Then $E_1$ is a reduction of $E_2$ in $F$ if and only if ${\rm e}_{\rm BR}(E_1)={\rm e}_{\rm BR}(E_2)$.
\end{lemma}

\subsection{Buchsbaum-Rim multiplicity in arbitrary degree}\label{section4.2} 
In analogy with the Subsection \ref{secBR}, in which we defined the classical Buchsbaum-Rim multiplicity, we consider the Buchsbaum-Rim multiplicity to a higher degree. Fixed some integer $e>0$. Using the notation described in Section \ref{section2}, let $(R,\mathfrak{m})$ be a $d$-dimensional Noetherian local ring, $E$ be an $R$-submodule of $F^e:=S_e$ such that $\ell_R(F^{e}/E) < \infty$. Let $N$ be an $R$-module of dimension $d$. For integers $n$ sufficiently large, Kleiman and Thorup in \cite[Section 8 and Proposition 8.2]{Kleiman-Thorup2}, showed the existence of an integer ${\rm e}^{d+p-1,0}(E)\in \mathbb{Z}$  such that $\ell_R\left(\frac{F^{en}\otimes_RN}{E^n\otimes_RN}\right)$ can be expressed as a polynomial in the form:
$$P_E^F\left(n;N\right)={\rm e}^{d+p-1,0}\left(E;N\right)\frac{n^{d+p-1}}{\left(d+p-1\right)!}+\mbox{lower terms}.$$
The integer ${\rm e}^{d+p-1,0}(E;N)$ is called the {\it Buchsbaum-Rim multiplicity in higher degree of $E$ in $F^e$}. In this paper, for simplicity, we will denote ${\rm e}^{d+p-1,0}(E;N)$ by ${\rm \tilde{e}}_{\rm BR}(E;N)$. Note that if $e=1$, the integer ${\rm e}^{d+p-1,0}(E;N)$ coincides with the usual Buchsbaum-Rim multiplicity ${\rm e}_{\rm BR}(E;N)$.

\subsection{Mixed multiplicity for modules} \label{section4.3} Let $(R,\mathfrak{m})$ be a $d$-dimensional Noetherian local ring and $E_1, \ldots, E_k$ be $R$-submodules of $F$ with $\ell_R(F/E_i) < \infty$, for $i = 1, \ldots, k$.  Let $N$ be an $R$-module of dimension $d$. The length function $h_R(n_1, \ldots, n_k; N) = \ell_R\left(\frac{F^{|{\bf n}|}\otimes_RN}{{\bf E}^{{\bf n}}\otimes_RN}\right)$ is a polynomial of total degree at most $d+p-1$, for sufficiently large values of $n_1, \ldots, n_k > 0$. The leading term of this polynomial can be expressed as:
$$\sum_{|{\bf d}|=d+p-1}\;\frac{1}{{\bf d}!}{\rm e}_{\rm BR}\left(E_1^{[d_1]},
\ldots, E_k^{[d_k]};N\right){\bf n}^{\bf d}.$$

Here, ${\bf d} = (d_1, \ldots, d_k)$ is a multi-index such that $|{\bf d}| = d+p-1$. The coefficients ${\rm e}_{\rm BR}\left(E_1^{[d_1]}, \ldots, E_k^{[d_k]};N\right)$ are called of {\it mixed Buchsbaum-Rim multiplicity of $E_1, \ldots, E_k$ of type $(d_1, \ldots, d_k)$ with respect to $N$}. In this definition, $E_i^{[d_i]}$ denotes that the module $E_i$ is listed $d_i$ times. This definition generalizes the notion of mixed multiplicities of $\mathfrak{m}$-primary ideals.

The definition of ${\rm e}_{\rm BR}\left(E_1^{[d_1]}, \ldots, E_k^{[d_k]};N\right)$ is consistent with the definitions given by Bedregal-Perez \cite[Section 4]{Bedregal-Perez2}, Kirby-Rees \cite{Kirby-Rees1}, and Kleiman-Thorup \cite[Section 8]{Kleiman-Thorup2}. In the special case, when $k = d+p-1$ and $d_1 = \cdots = d_k = 1$, we can denote ${\rm e}_{\rm BR}\left(E_1^{[d_1]}, \ldots, E_k^{[d_k]};N\right)$ as ${\rm e}_{\rm BR}\left(E_1, \ldots, E_{d+p-1};N\right)$.

If $N = R$, we use the notation ${\rm e}_{\rm BR}\left(E_1^{[d_1]}, \ldots, E_k^{[d_k]}\right)$ to denote ${\rm e}_{\rm BR}\left(E_1^{[d_1]}, \ldots, E_k^{[d_k]};R\right)$.

\subsection{Mixed multiplicity for modules for arbitrary degree} \label{MixHdegree} Fixed some integers $e_i\in \mathbb{N}$. Let $(R,\mathfrak{m})$ be a $d$-dimensional Noetherian local ring  and $E_1, \ldots, E_k$ be $R$-submodules of $F^{e_1}, \ldots, F^{e_k}$ such that $\ell_R(F^{e_i}/E_i) < \infty$ for all $i = 1, \ldots, k$. Let $N$ be an $R$-module of dimension $d$. Then, the length function $h_R\left(n_1, \ldots, n_k; N\right) = \ell_R\left(\frac{F^{{\bf e\cdot n}}\otimes_RN}{{\bf E}^{{\bf n}}\otimes_RN}\right)$, for sufficiently large values of $n_1, \ldots, n_k > 0$, can be expressed as a polynomial of total degree at most $d+p-1$, and its leading term is given by:
$$\sum_{|{\bf d}|=d+p-1}\;\frac{1}{{\bf d}!}{\rm \tilde{e}}_{\rm BR}\left(E_1^{[d_1]},
\ldots, E_k^{[d_k]};N\right){\bf n}^{\bf d},$$
where ${\bf d} = (d_1, \ldots, d_k)$ is a multi-index with $|{\bf d}| = d+p-1$. The coefficients ${\rm \tilde{e}}_{\rm BR}(E_1^{[d_1]}, \ldots, E_k^{[d_k]};N)$ are called {\it mixed Buchsbaum-Rim multiplicity of $E_1, \ldots, E_k$ of type $(d_1, \ldots, d_k)$ with respect to $N$}. This definition agrees with the one given in \cite[p. 568]{Kleiman-Thorup2}.

Moreover, it is important to note that 
 when $e_1 = e_2 = \cdots = e_k = 1$, the mixed Buchsbaum-Rim multiplicity is reduced to the usual mixed multiplicity, that is, 
$${\rm \tilde{e}}_{\rm BR}\left(E_1^{[d_1]}, \ldots, E_k^{[d_k]}; N\right) = {\rm e}_{\rm BR}\left(E_1^{[d_1]}, \ldots, E_k^{[d_k]};N\right).$$

 \subsection{Multiplicity system} 
 In this subsection, we will define the $g$-multiplicity system introduced by Kirby in \cite[Proposition 2]{Kirby}. In our context, this pertains to the ring $S = \operatorname{Sym}_R(F)$. To do this, we will recall the definition of the Koszul complex $K_{\bullet}(a_1, \ldots, a_m; N)$ of a set of homogeneous elements $a_1, \ldots, a_m$ in $S$ with $\deg(a_i) = k_i$ for $i = 1, \ldots, m$, with respect to an $R$-module $N.$

  Let $G = \bigoplus_{i=1}^m S(-k_i)$ be a free $S$-module with basis $e_1, \dots, e_m$. The homomorphism $\phi: G \to S$ defined by $\phi(e_i) = a_i$ gives rise to the Koszul complex $K_{\bullet}(a_1,\dots, a_m;S)$, where the $n$-th graded piece is $K_n(a_1,\dots, a_m;S) = \bigwedge^n G$ and the differential $d_n: \bigwedge^n G \to \bigwedge^{n-1} G$ is defined as
$$d_n(x_1 \wedge \cdots \wedge x_n) = \sum_{i=1}^n (-1)^{i+1} \phi(x_i)x_1 \wedge \cdots \wedge \widehat{x_i} \wedge \cdots \wedge x_n,$$
where $\widehat{x_i}$ denotes the omission of $x_i$.

By identifying the basis elements $e_{i_1} \wedge \cdots \wedge e_{i_n}$ with $e_{{i_1,\dots,i_n}}$, the complex $K_{\bullet}(a_1,\dots, a_m;S)$ can be written as $K_n(a_1,\dots, a_m;S) = \bigoplus_{1\leq i_1<\cdots <i_n\leq m} Se_{{i_1,\dots,i_n}}$, where $\deg(e_{{i_1,\dots,i_n}}) = k_{i_1} + \dots + k_{i_n}$. This complex has a graded structure, and the differential $d_n$ is an homogeneous morphism of graded modules.

Let $N$ be an $R$-module of dimension $d$. For the graded $S$-module $M = S \otimes_R N = \oplus_{q\geq 0} F^q \otimes_R N$, the homological Koszul complex $K_{\bullet}(a_1,\dots, a_m;N)$ is obtained by tensoring $K_{\bullet}(a_1,\dots, a_m;S)$ with $N$, where $d_n \otimes \operatorname{id}_N$ defines the differentials. The complex $K_{\bullet}(a_1,\dots, a_m;N)$ is a graded $S$-module, with each component $[K_{\bullet}(a_1,\dots, a_m;N)]_t$ corresponding to a complex of $R$-modules.

The Koszul homology modules $H_n(a_1,\dots, a_m;N)$ are defined as the homology modules of the Koszul complex, that is, 
$$H_n\left(a_1,\dots, a_m;N\right) = \frac{\operatorname{Ker}\left(d_n \otimes \operatorname{id}_N\right)}{\operatorname{Im}\left(d_{n+1} \otimes \operatorname{id}_N\right)},$$ 
for $n \geq 0$. These modules are finitely generated graded $S$-modules.

The $i$-th homology module $H_iK_t(a_1,\dots, a_m;N)$ of the component $[K_{\bullet}(a_1,\dots, a_m;N)]_t$ is denoted by $H_iK_t(a_1,\dots, a_m;N)$.

When the set of homogeneous elements $a_1,\dots,a_m$ of $S$ form a {\it $g$-multiplicity system} with respect to $N$, the components $(M/\sum_{i=1}^m a_i M)_t$ have finite $R$-length for sufficiently large $t\geq 0$. In this case, the Koszul homology modules $H_iK_t(a_1,\dots, a_m;N)$ also have finite length for sufficiently large $t\geq 0$. D. Kirby \cite[Proposition 2]{Kirby} introduced the {\it $g$-multiplicity} of degree $t$ with respect to $N$, denoted by ${\rm e_t}(a_1,\dots,a_m;N)$, defined as
$${\rm e_t}((a_1,\dots,a_m);N) = \sum_{i=0}^m (-1)^i \ell_R(H_iK_t(a_1,\dots,a_m ;N)).$$

\subsection{ Some properties about multiplicities} 
In this subsection, we state results regarding multiplicities and provide properties that establish connections between previously defined multiplicities. These results carry significant importance in the subsequent sections, where they will serve as fundamental tools in proving the main result. One important result in this context is the associativity formula for the Buchsbaum-Rim multiplicity, which has been proved by Kleiman  \cite[Proposition 7]{Kleiman2}. Similarly, a formula for the Buchsbaum-Rim multiplicity can be found in \cite[Corollary 2.8]{Bedregal-Perez2} and \cite[Theorem 4.5(v)]{Kirby-Rees1}. Furthermore, Validashti \cite[Theorem 6.5.1]{Validashti}, establishes an associativity formula for the $j$-multiplicity, which generalizes the Buchsbaum-Rim multiplicity.

These associativity formulas provide a way to compute multiplicities or mixed multiplicities on domains. Specifically, they are applicable to reduced Noetherian local rings that are domains. 

\begin{lemma}[Associativity formula for Buchsbaum-Rim multiplicity] \label{Associativity} Let $(R,\mathfrak{m})$ be a $d$-dimensional Noetherian local ring, $E$ be an $R$-submodule of $F$ such that $\ell_R(F/E) <\infty$ and $N$ be an $R$-module of dimension $d$. Define $\Lambda_d=\{\mathfrak{p}\in \Supp(N) : \dim( R/\mathfrak{p})=d\}$. Then 
$${\rm e}_{\rm BR}\left(E;N\right)=\sum_{\mathfrak{p}\in \Lambda_d }\ell_{R_{\mathfrak{p}}}\left(N_{\mathfrak{p}}\right)
{\rm e}_{\rm BR}\left(E\left(\frac{R}{\mathfrak{p}}\right);\frac{R}{\mathfrak{p}}\right).$$
\end{lemma}

The associativity formula for the mixed Buchsbaum-Rim multiplicity is also well-known, see Bedregal-Perez \cite[Proposition 5.4]{Bedregal-Perez2} or Kirby-Rees \cite[Theorem 6.3(iii)]{Kirby-Rees1}.

\begin{lemma}[Associativity formula for mixed multiplicities] \label{AssociativityMix} Let $(R,\mathfrak{m})$ be a $d$-dimensional Noetherian local ring, $E_1, \ldots, E_k$ be $R$-submodules of $F$ such that $\ell_R(F/E_i) < \infty$ for all $i = 1, \ldots, k$, and let $N$ be an $R$-module of dimension $d$. Let $\Lambda_d = \{\mathfrak{p}\in \Supp(N) : \dim(R/\mathfrak{p}) = d\}$. Then for any integers $d_1, \ldots, d_k$ with $|\mathbf{d}| = d + p - 1$,
$${\rm e}_{\rm BR}\left(E^{[d_1]}_1, \ldots, E^{[d_k]}_k; N\right) = \sum_{\mathfrak{p} \in \Lambda_d} \ell_{R_{\mathfrak{p}}}\left(N_{\mathfrak{p}}\right) {\rm e}_{\rm BR} \left(E^{[d_1]}_1\left(\frac{R}{\mathfrak{p}}\right), \ldots, E^{[d_k]}_k\left(\frac{R}{\mathfrak{p}}\right); \frac{R}{\mathfrak{p}}\right).$$
\end{lemma}

\begin{lemma}\label{Lemma 17.5.3} 
Let $(R,\mathfrak{m})$ be a $d$-dimensional Noetherian local ring, $E_1, \dots, E_{d+p-1}$ be $R$-submodules of $F$ such that $\ell_R(F/E_i) < \infty$ for all $i=1, \dots, d+p-1$, and let $N$ be an $R$-module of dimension $d$.
\begin{itemize}
\item[(i)] If $L_1, \dots, L_{d+p-1}$ are $R$-submodules of $F$ with $\ell_R(F/L_i) < \infty$ for $i=1, \dots, d+p-1$, and $L_i \subseteq E_i$. Then
$${\rm e}_{\rm BR}\left(L_1, \ldots, L_{d+p-1}; N\right) \geq {\rm e}_{\rm BR}\left(E_1, \ldots, E_{d+p-1}; N\right).$$

\item[(ii)] If $x_i \in E_i$, for $i = 1, \dots, d+p-1$, and $(x_1, \dots, x_{d+p-1})$ is an $R$-module with $\ell_R(F/(x_1, \dots, x_{d+p-1})) < \infty$. Then
$${\rm e}_{\rm BR}\left(\left(x_1, \dots, x_{d+p-1}\right); N\right) \geq {\rm e}_{\rm BR}\left(E_1, \dots, E_{d+p-1}; N\right).$$
\end{itemize}
\end{lemma}

\begin{proof}
The proof of (i) follows from \cite[Corollary 6.4]{Bedregal-Perez2}. To prove (ii), first apply Proposition \ref{Proposition 17.3.3} to construct $R$-modules of finite length, $L_1, \dots, L_{d+p-1}$, with $x_i \in L_i \subseteq E_i$ for all $i$, and $(x_1, \dots, x_{d+p-1})$ being a joint reduction of $L_1, \dots, L_{d+p-1}$ with respect to $N$. Then (ii) follows from (i) and \cite[Theorem 5.5]{Bedregal-Perez2}.
\end{proof}

\begin{proposition}\label{Proposition 11.2.9}
Let $(R,\mathfrak{m})$ be a $d$-dimensional Noetherian local ring, $E$ be an $R$-submodule of $F$ such that $\ell_R(F/E)<\infty$ and let $N$ be an $R$-module of dimension $d$.
\begin{enumerate}
\item[(i)] If $r$ is a positive integer, then ${\rm \tilde{e}}_{\rm BR}(E^r) = {\rm e}_{\rm BR}(E)r^{d+p-1}$.
\item[(ii)] If $E = (x_1, \dots, x_{d+p-1})$, then for any $l_1, \dots, l_{d+p-1} \in \mathbb{N}$, $${\rm e}_{t}\left(x_1^{l_1}, \dots, x_{d+p-1}^{l_{d+p-1}}; N\right) = l_1 \cdots l_{d+p-1} \cdot {\rm e_{BR}}\left((x_1, \dots, x_{d+p-1}); N\right),$$
for sufficiently large $t>0$.
\end{enumerate}
\end{proposition}

\begin{proof}
(i) For $n$ sufficiently large, 
$$\ell_R\left(\frac{F^{rn}}{\left(E^r\right)^n}\right) = {\rm \tilde{e}}_{\rm BR}(E^r)\frac{n^{d+p-1}}{\left(d+p-1\right)!} + \text{lower terms}.$$

On the other hand, for $n$ sufficiently large, 
$$\begin{array}{rcl}
\displaystyle \ell_R\left(\frac{F^{rn}}{E^{rn}}\right) 
& = &  \displaystyle {\rm e}_{\rm BR}(E)\frac{\left(rn\right)^{d+p-1}}{\left(d+p-1\right)!} + \text{lower terms} \\
\\
& = & \displaystyle {\rm e}_{\rm BR}(E)r^{d+p-1}\frac{n^{d+p-1}}{\left(d+p-1\right)!} + \text{lower terms}.
\end{array}$$

Now, comparing the coefficients of $n^{d+p-1}$, we obtain the desired equality.

(ii) By \cite[Lemma 2.7(i)]{Bedregal-Perez3}, 
$${\rm e}_{t}\left(x_1^{l_1}, \dots, x_{d+p-1}^{l_{d+p-1}}; N\right) = l_1 \cdots l_{d+p-1} \cdot {\rm e_{t}}\left(x_1, \dots, x_{d+p-1}; N\right).$$

Since $x_i$ are homogeneous elements of degree one, by \cite[Theorem 4.1]{Bedregal-Perez3}, follows ${\rm e_{t}}(x_1, \dots, x_{d+p-1}; N) = {\rm e_{\rm BR}}((x_1, \dots, x_{d+p-1}); N)$.
\end{proof}

\section{The Risler-Teissier Mixed Multiplicity Theorem}\label{section4}

In this section, we state and prove the {\it Risler-Teissier Theorem} for finitely generated $R$-submodules $E_i$ of $F^{e_i}$, when $e_i>0$ are positive integers. This result generalizes the Risler-Teissier Theorem for ideals given in \cite[Proposition 2.1]{Teissier}, and also generalizes the Risler-Teissier theorem for modules given in \cite[Theorem 5.2]{Bedregal-Perez2} when $e_i=1$, for all $i$. For this, let us first recall the definition of associated mixed Buchsbaum-Rim multiplicity introduced by Kirby-Rees in \cite{Kirby-Rees1}  and by Kleiman-Thorup in \cite[p. 532 and  568]{Kleiman-Thorup2}. 

Let $(R,\mathfrak{m})$ be a $d$-dimensional Noetherian local ring, $E_1, \ldots, E_k$ be $R$-submodules of $F^{e_1}, \ldots, F^{e_k}$ such that $\ell_R(F^{e_i}/E_i) < \infty$ for all $i = 1, \ldots, k$, $e_i \in \mathbb{N}$ be fixed integers. Let $N$ be an $R$-module of dimension $d$, and consider the graded $S$-module $M=S\otimes_R N$. The length function
$$h_R(n_1, \ldots, n_k,q; N)=\ell_R\left(\frac{M_{{\bf e \cdot n}+q}}{{\bf E}^{{\bf n}}{M_q}}\right),$$
for sufficiently large $n_1, \ldots , n_k, q$, is a polynomial of total degree at most $d+p-1$, whose leading term can be written as
$$\sum_{j+|{\bf d}|=d+p-1}\;\frac{1}{{\bf d}!}{\rm \tilde{e}}_{\rm BR}^j\left(E_1^{[d_1]},
\ldots, E_k^{[d_k]};N\right){\bf n}^{\bf d}q^j,$$
where ${\bf d} = (d_1, \ldots, d_k)$ is a multi-index with $j+|{\bf d}| = d+p-1$. The coefficients ${\rm \tilde{e}}_{\rm BR}^j\left(E_1^{[d_1]}, \ldots, E_k^{[d_k]};N\right)$ are called  the {\it associated mixed Buchsbaum-Rim multiplicity} of $E_1, \ldots, E_k$ of type $(d_1, \ldots, d_k)$ with respect to $N$. Note that, when $j=0$, the terms ${\rm \tilde{e}}_{\rm BR}^0\left(E_1^{[d_1]}, \ldots, E_k^{[d_k]};N\right)$ are the mixed Buchsbaum-Rim multiplicity 
${\rm \tilde{e}}_{\rm BR}\left(E_1^{[d_1]}, \dots, E_k^{[d_k]};N\right)$ as defined in Subsection \ref{MixHdegree} or see also \cite[p. 568]{Kleiman-Thorup2}.

\begin{lemma} \label{Lemma 17.2.4Degree} Let $R$ be a Noetherian ring and let $N$ be an $R$-module. Let $E_1,\dots,E_k$ be $R$-submodules of $F^{e_i}$ with positive integers $e_i$, for all $i=1,\dots,k$. Suppose that $x \in E_1$ is a superficial element for $E_1, E_2, \ldots, E_k$ with respect to $N$. Assume that $I_{E_1}\subseteq  \sqrt{I_{E_2}\cdots I_{E_k}}$ and $\cap_n I_{E_1}^n M = 0$, with $M=S\otimes_R N$. Then, for positive integer $q>0$,
{\small{$$\left({\bf E}^{{\bf n}}M_q :_{M_{{\bf e\cdot n}-e_1+q}} x\right) = (0 :_{M_{{\bf e\cdot n}-e_1+q}} x) + {\bf E}^{{\bf n}-e_1\delta(1)}M_q\,\,\ \mbox{  and }\,\,\,\left(0 :_{M_{{\bf e \cdot n}-e_1+q}} x\right)\cap {\bf E}^{{\bf n}-e_1\delta(1)}M_q=0.$$}}
\end{lemma}
\begin{proof} To prove this lemma, it is sufficient to consider the ideals $I_{E_1},\dots,I_{E_k}$ of $S$, for all $i=1,\dots,k$. Thus, according to \cite[Lemma 17.2.4]{SH}, for all sufficiently large $n_1,\dots,n_k$, we obtain

{\small{$$({\bf I}_{\bf E}^{{\bf n}}M :_M x) = (0 :_M x) + {\bf I}_{\bf E}^{{\bf n}-e_1\delta(1)}M\,\,\ \mbox{ and }\,\,\,\left(0 :_{M} x\right)\cap {\bf I}_{\bf E}^{{\bf n}-e_1\delta(1)}M=0.$$}}

Hence, the result follows by Remark \ref{relasuperficial} and by concentrating on degree ${\bf e\cdot n}-e_1+q$ in the above equalities.

\end{proof}

\begin{theorem}[The Risler-Teissier mixed multiplicity theorem for $e_i\geq 1$]\label{Theorem 17.4.6degree} Let $(R,\mathfrak{m})$ be a $d$-dimensional Noetherian local ring  with infinite residue field and let $N$ be an $R$-module of dimension $d$. Let $E_1,\dots, E_{k}$ be $R$-submodules of $F^{e_1},\dots, F^{e_k}$, respectively, such that, $\ell_R(F^{e_i}/E_i)<\infty$, $e_i\in \mathbb{N}$ for all $i=1,\dots,k$. Fix $c_1,\ldots,c_k\in\mathbb{N}$ and assume $e_1=\min\{e_1,\dots,e_k\}$. 
Let $x \in E_1$ be a superficial element for $E_1,\ldots,E_k$ with respect to $N$ and denote the $S$-graded module $M=S\otimes_RN$. Assume that $x$ is not contained in any minimal prime of ${\rm Ann}_S(M)$. 
Set $S'= S/xS$ and $M'= M/xM$. Then, for any integers $j, d_1,\ldots,d_k\in\mathbb{N}$ 
with $j+\vert{\bf d}\vert=d+p-1$, $d_1>0$, and for sufficiently large integers $n_1,\dots,n_k,$ $q$,
{\small{$$
{\rm \tilde{e}}_{\rm BR}^j\left(E_1^{[d_1]},\ldots,E_k^{[d_k]};N\right)=
\left\{\begin{array}{ll}
\frac{1}{e_1}{\rm \tilde{e}}_{\rm BR}^j\left({E_1'}^{[d_1-1]},\ldots,{E_k'}^{[d_k]};M'\right), & \textrm{if}\ d+p > 2; \\ \\
\frac{1}{e_1}\left[\ell_R\left(\frac{F^{{\bf e\cdot n}+q}\otimes_RN}{(x)F^{{\bf e\cdot n}-e_1+q}\otimes_RN}\right) - \ell_R\left(0:_{F^{{\bf e\cdot n}-e_1+q}\otimes_RN} x\right)\right], & \textrm{if}\ d+p = 2.
\end{array}
\right.
$$
}}
where $E'_i$ is image of  $E_i$ in  $F^{e_i}/(x)F^{e_i-e_1}$ for all  $i=1,\dots, k$. 
\end{theorem}
\begin{proof}
By the exact sequence
{\small{$$
 0 
 \longrightarrow \frac{\left({\bf E}^{{\bf n}}M_q:_{M_{{\bf e \cdot n}-e_1+q}} x\right)}{{\bf E}^{{\bf n}-e_1\delta(1)}{M_q}} \stackrel{i}{\longrightarrow} \frac{M_{{\bf e\cdot n}-e_1+q}}{{\bf E}^{{\bf n}-e_1\delta(1)}M_q}
 \stackrel{\cdot x}{\longrightarrow} \frac{M_{{\bf e \cdot n}+q}}{{\bf E}^{{\bf n}}{M_q}} \stackrel{\pi}{\longrightarrow} \frac{M_{{\bf e \cdot n}+q}}{xM_{{\bf e \cdot n}-e_1+q}+{\bf E}^{{\bf n}}{M_q}} \longrightarrow 0,
$$}}
\noindent we obtain

{\small{$$\displaystyle \ell_R\left(\frac{M_{{\bf e \cdot n}+q}}{{\bf E}^{{\bf n}}{M_q}}\right)-\ell_R\left(\frac{M_{{\bf e \cdot n}-e_1+q}}{{\bf E}^{{\bf n}-e_1\delta(1)}{M_q}}\right)=\ell_R\left(\frac{M_{{\bf e \cdot n}+q}}{xM_{{\bf e \cdot n}-e_1+q}+{\bf E}^{{\bf n}}M_q}\right)-\ell_R\left(\frac{({\bf E}^{{\bf n}}M_q :_{M_{{\bf e \cdot n}-e_1+q}} x)}{{\bf E}^{{\bf n}-e_1\delta(1)}M_q}\right).$$}}

Thus, by Lemma \ref{Lemma 17.2.4Degree}, for large $n_1,\dots ,n_k$, say $n_1 \geq c_1, \dots ,n_k \geq c_k$, and all $q\geq 0,$ we get $\left(0 :_{M_{{\bf e \cdot n}-e_1+q}} x\right)\cong \frac{({\bf E}^{{\bf n}}M_q :_{M_{{\bf e\cdot n}-e_1+q}} x)}{{\bf E}^{{\bf n}-e_1\delta(1)}{M_q}}.$ Then
{\small{\begin{equation}\label{Equa21}
\displaystyle \ell_R\left(\frac{M_{{\bf e \cdot n}+q}}{{\bf E}^{{\bf n}}{M_q}}\right)-\ell_R\left(\frac{M_{{\bf e \cdot n}-e_1+q}}{{\bf E}^{{\bf n}-e_1\delta(1)}{M_q}}\right)=\ell_R\left(\frac{M_{{\bf e \cdot n}+q}}{xM_{{\bf e \cdot n}-e_1+q}+{\bf E}^{{\bf n}}M_q}\right)-\ell_R\left(0 :_{M_{{\bf e \cdot n}-e_1+q}} x\right).
\end{equation}}}

Then if $d+p=2$, we have that $M'_{{\bf e\cdot n}+q}=\frac{M_{{\bf e \cdot n}+q}}{xM_{{\bf e \cdot n}-e_1+q}+{\bf E}^{{\bf n}}M_q}$ and the left side of Equation (\ref{Equa21}) has coefficient ${\rm \tilde{e}}_{\rm BR}^j\left(E_1^{[d_1]},\ldots,E_k^{[d_k]};N\right)$ in its maximal dimension.

Suppose now that $d+p > 2$ and consider the following exact sequences

\begin{equation}\label{Equa23}  
0\longrightarrow \frac{{\bf E}^{{\bf n}+e_1\delta(1)}M_{q-e_1}}{x{\bf E}^{{\bf n}}{M_{q-e_1}}} \stackrel{i}{\longrightarrow} \frac{{\bf E}^{{\bf n}}M_{q}}{x{\bf E}^{{\bf n}}{M_{q-e_1}}} \stackrel{\pi}{\longrightarrow} \frac{{\bf E}^{{\bf n}}M_{q}}{{\bf E}^{{\bf n}+e_1\delta(1)}{M_{q-e_1}}}\longrightarrow 0
\end{equation}
and
\begin{flushleft} 
$\displaystyle 0\longrightarrow \frac{x\left({\bf E}^{{\bf n}+e_1\delta(1)}M_{q-e_1}:_{{\bf E}^{{\bf n}-e_1\delta(1)}M_q}x\right)}{x{\bf E}^{{\bf n}}{M_{q-e_1}}} \longrightarrow \frac{{\bf E}^{{\bf n}+e_1\delta(1)}M_{q-e_1}}{x{\bf E}^{{\bf n}}{M_{q-e_1}}} \stackrel{i}{\longrightarrow}$
\end{flushleft}
\begin{flushright}
$\displaystyle \frac{{\bf E}^{{\bf n}}M_{q}}{x{\bf E}^{{\bf n}-e_1\delta(1)}{M_{q}}} \stackrel{\pi}{\longrightarrow} \frac{{\bf E}^{{\bf n}}M_{q}}{x{\bf E}^{{\bf n}-e_1\delta(1)}{M_{q}}+{\bf E}^{{\bf n}+e_1\delta(1)}{M_{q-e_1}}}\longrightarrow 0.$
\end{flushright}

Consider the homomorphism:
$$\phi:\frac{\left({\bf E}^{{\bf n}+e_1\delta(1)}M_{q-e_1}:_{{\bf E}^{{\bf n}-e_1\delta(1)}M_q}x\right)}{{\bf E}^{{\bf n}-e_1\delta(1)}{M_{q}}}\longrightarrow \frac{x\left({\bf E}^{{\bf n}+e_1\delta(1)}M_{q-e_1}:_{{\bf E}^{{\bf n}-e_1\delta(1)}M_q}x\right)}{x{\bf E}^{{\bf n}-e_1\delta(1)}{M_{q}}}$$
with $\ker(\phi)=\{0\}$ since $x \in E_1$ is a superficial element and Lemma \ref{Lemma 17.2.4Degree}. 
Then $\frac{\left({\bf E}^{{\bf n}+e_1\delta(1)}M_{q-e_1}:_{{\bf E}^{{\bf n}-e_1\delta(1)}{M_{q}}}x\right)}{{\bf E}^{{\bf n}-e_1\delta(1)}{M_{q}}} \cong 
\frac{x\left({\bf E}^{{\bf n}+e_1\delta(1)}{M_{q-e_1}}:_{{\bf E}^{{\bf n}-e_1\delta(1)}{M_{q}}}x\right)}{x{\bf E}^{{\bf n}-e_1\delta(1)}{M_{q}}}$.

Similarly, we also get that $\frac{{\bf E}^{{\bf n}-e_1\delta(1)}{M_{q}}}{{\bf E}^{{\bf n}}{M_{q-e_1}}}\cong\frac{x{\bf E}^{{\bf n}-e_1\delta(1)}{M_{q}}}{x{\bf E}^{{\bf n}}{M_{q-e_1}}}$.

We obtain the following exact sequence:
{\tiny{
\begin{equation}\label{Equa27} \xymatrix{ 
  & \frac{\left({\bf E}^{{\bf n}+e_1\delta(1)}M_{q-e_1}:_{{\bf E}^{{\bf n}-e_1\delta(1)}{M_{q}}}x\right)}{{\bf E}^{{\bf n}-e_1\delta(1)}{M_{q}}}  \ar[d]^{\phi}&  & &  & \\
0 \ar[r] & \frac{x\left({\bf E}^{{\bf n}+e_1\delta(1)}{M_{q-e_1}}:_{{\bf E}^{{\bf n}-e_1\delta(1)}{M_{q}}}x\right)}{x{\bf E}^{{\bf n}-e_1\delta(1)}{M_{q}}}\ar[r] & \frac{{\bf E}^{{\bf n}+e_1\delta(1)}M_{q-e_1}}{x{\bf E}^{{\bf n}-e_1\delta(1)}{M_{q}}} \ar[r] & \frac{{\bf E}^{{\bf n}}M_{q}}{x{\bf E}^{{\bf n}-e_1\delta(1)}{M_{q}}}  \ar[r] & \frac{{\bf E}^{{\bf n}}M_{q}}{x{\bf E}^{{\bf n}-e_1\delta(1)}{M_{q}}+{\bf E}^{{\bf n}+e_1\delta(1)}{M_{q-e_1}}} \ar[r] & 0
} 
\end{equation}}}
and 
\begin{equation}\label{Equa26}
    0 \longrightarrow \frac{{\bf E}^{{\bf n}-e_1\delta(1)}{M_{q}}}{{\bf E}^{{\bf n}}{M_{q-e_1}}} \longrightarrow \frac{{\bf E}^{{\bf n}-e_1\delta(1)}{M_{q}}}{x{\bf E}^{{\bf n}}{M_{q-e_1}}} \longrightarrow 
\frac{{\bf E}^{{\bf n}-e_1\delta(1)}{M_{q}}}{x{\bf E}^{{\bf n}-e_1\delta(1)}{M_{q}}} \longrightarrow 0.
\end{equation}

Therefore from the exact sequences (\ref{Equa23}), (\ref{Equa27}) and (\ref{Equa26}) we have that

\begin{equation}\label{Equ29}
    \ell_R\left(\frac{{\bf E}^{{\bf n}}M_{q}}{x{\bf E}^{{\bf n}-e_1\delta(1)}{M_{q}}+{\bf E}^{{\bf n}+e_1\delta(1)}{M_{q-e_1}}}\right)=\ell_R\left(\frac{{\bf E}^{{\bf n}}M_{q}}{{\bf E}^{{\bf n}+e_1\delta(1)}{M_{q-e_1}}}\right)-\ell_R\left(\frac{{\bf E}^{{\bf n}-e_1\delta(1)}{M_{q}}}{{\bf E}^{{\bf n}}{M_{q-e_1}}}\right).
\end{equation}

By the definition above, for sufficiently large $n_1, \dots, n_k, q$, each term on the left side of (\ref{Equ29}), for $j=0$, is a polynomial whose leading term has the form:

$$\sum_{|{\bf d}|=d+p-1}\;\frac{1}{{\bf d}!}{\rm \tilde{e}}_{\rm BR}^0\left(E_1^{[d_1]},
\ldots, E_k^{[d_k]};N\right)\left({n_1}^{d_1}-{\left(n_1+e_1\right)}^{d_1}\right){n_2}^{d_2}\dots{n_k}^{d_k}$$
and
$$\sum_{|{\bf d}|=d+p-1}\;\frac{1}{{\bf d}!}{\rm \tilde{e}}_{\rm BR}^0\left(E_1^{[d_1]},
\ldots, E_k^{[d_k]};N\right)\left({n_1}^{d_1}-{\left(n_1-e_1\right)}^{d_1}\right){n_2}^{d_2}\dots{n_k}^{d_k}.$$

Thus, the difference on the left side is:

$$\sum_{|{\bf d}|=d+p-1}\;\frac{1}{(d_1-2)!d_2!\dots d_k!}{\rm \tilde{e}}_{\rm BR}^0\left(E_1^{[d_1]},
\ldots, E_k^{[d_k]};N\right)e_1n_1^{d_1-2}{n_2}^{d_2}\dots{n_k}^{d_k}.$$

Therefore, the proposition follows by comparing the leading coefficients.
\end{proof}

It should be noted that the proof of the previous theorem is similar, except for the submodules contained in free modules of degree $e_1\geq 1$, to the proof given in \cite[Theorem 5.2]{Bedregal-Perez2}. For this reason, we omit the demonstration of certain isomorphisms.

\begin{remark}\label{CorteRemark} Based on the definition of mixed multiplicities, we have the following relationship: if $k = 1$, then ${\rm \tilde{e}}_{\rm BR}(E_1^{[d_1]};N) = {\rm \tilde{e}}_{\rm BR}(E_1;N)$. Consequently, using the Risler-Teissier mixed multiplicity theorem \ref{Theorem 17.4.6degree}, we obtain:

$$
{\rm \tilde{e}}_{\rm BR}\left(E_1;N\right) =
\left\{
\begin{array}{ll}
\frac{1}{e_1}{\rm \tilde{e}}_{\rm BR}\left({E_1'};M'\right), & \textrm{if } d+p > 2, \\ \\
\frac{\ell_R\left(\frac{F^{e_1\cdot n}\otimes_RN}{(x)F^{e_1\cdot n-e_1}\otimes_RN}\right) - \ell_R\left(0:_{F^{e_1\cdot n-e_1}\otimes_RN} x\right)}{e_1}, & \textrm{if } d +p = 2.
\end{array}
\right.
$$

Here $E_1'=(E_1+(x))/(x)$ and $M'=S\otimes_RN/xS\otimes_RN$.
    
\end{remark}

\begin{corollary}\label{Theorem 17.4.6degreecor} Let $(R,\mathfrak{m})$ be a $d$-dimensional Noetherian local ring  with infinite residue field. Let $E_1,\dots, E_{k}$ be $R$-submodules of $F^e$ such that $\ell_R(F^e/E_i)<\infty$, for $l\in\mathbb{N}$ and all $i=1, \dots, k$. Let $N$ be a finitely generated $R$-module of dimension $d$. Let $x_1,\dots ,x_{d+p-1}$ be any superficial sequence for $E_1,\dots, E_{k},$ with respect to $N$, with each $E_i$ is listed $d_i$ times and each $(x_i)S$ is not in any minimal prime ideal over $(x_1 ,\dots , x_{i-1})S$. Then, for any integers $d_1,\ldots,d_k\in\mathbb{N}$ 
with $d_1+d_2+\cdots+d_k=d+p-1$, $d_i>0$, and for large $n$,
{\small{$$
{\rm \tilde{e}}_{\rm BR}\left(E_1^{[d_1]},\ldots,E_k^{[d_k]};N\right)=\frac{
\ell_R\left(\frac{M_{en}}{\left(x_1,\dots,x_{d+p-1}\right)M_{en-e}}\right) - \ell_R\left((x_1,\dots,x_{d+p-1})M_{en}:_{M_{en}} x_{d+p-1}\right)}{e^{d+p-1}},
$$}}
where 
$M_n=F^n\otimes_RN$, which is equal to $
{\rm\tilde{e}_{BR}}(x_1 ,\dots , x_{d+p-1} ; N )$.
\end{corollary}
\begin{proof}
To prove the first statement, we can use Theorem \ref{Theorem 17.4.6degree} iteratively.

Now, we proof the second statement. Suppose $d+p-1 > 1$, and let $E = (x_1,\dots,x_{d+p-1})$.  By repeatedly utilizing Remark \ref{CorteRemark}, we get 
{\small{$$
{\rm \tilde{e}}_{\rm BR}\left(E;N\right)=\frac{
\ell_R\left(\frac{M_{en}}{\left(x_1,\dots,x_{d+p-1}\right)M_{en-e}}\right) - \ell_R\left((x_1,\dots,x_{d+p-1})M_{en}:_{M_{en}} x_{d+p-1}\right)}{e^{d+p-1}},
$$}}
thus, it proves the corollary.
\end{proof}

In particular, if $e=1$, from Corollary \ref{Theorem 17.4.6degreecor}  we get the same result given in \cite[Theorem 5.2, Corollary 5.3]{Bedregal-Perez2}.

\begin{corollary}\label{Theorem 17.4.6degreecor1} Let $(R,\mathfrak{m})$ be a $d$-dimensional Noetherian local ring  with infinite residue field. Let $E_1,\dots, E_{k},$ be $R$-submodules of  $F$  such that  $F/E_i$ has  finite length  and $e\in \mathbb{N}$, for all $i=1, \dots, k$. Let $N$ be a finitely generated $R$-module of dimension $d$. Let $x_1,\dots ,x_{d+p-1}$ be any superficial sequence for $E_1,\dots, E_{k},$ with respect to $N$ , with each $E_i$ listed $d_i$ times and each $(x_i)S$ is not in any minimal prime ideal over $(x_1 ,\dots , x_{i-1})S$. Then, for any integers $d_1,\ldots,d_k\in\mathbb{N}$ 
with $d_1+d_2+\cdots+d_k=d+p-1$, $d_i>0$, and for large $n$,
{\small{$$
{\rm e}_{\rm BR}\left(E_1^{[d_1]},\ldots,E_k^{[d_k]};N\right)=
\ell_R\left(\frac{M_{n}}{(x_1,\dots,x_{d+p-1})M_{n-1}}\right) - \ell_R\left(\left(x_1,\dots,x_{d+p-1}\right)M_{n}:_{M_n} x_{d+p-1}\right),
$$}}
where 
$M_n=F^n\otimes_RN$, which is equal to ${\rm e}_{\rm BR}(x_1 ,\dots , x_{d+p-1} ; N )$.
\end{corollary}

\section{Fundamental lemmas}\label{section5}

In this section, we present some fundamental results that will be crucial to prove the main theorem of this paper.

\begin{lemma}\label{Lemma 17.5.2} Let $(R, \mathfrak{m})$ be a Noetherian local ring with an infinite residue field. Let $E_1, \dots, E_k$ be finitely generated $R$-submodules of a free module $F$ with a positive rank $p$. Suppose $I_{E_1} \subseteq \sqrt{I_{E_2} \cdots I_{E_k}}$, and let $N$ be a finitely generated $R$-module. Consider a variable $Y$ over $R$ and an element $x \in E_1$. Then there exist positive integers $c$ and $e$, and a non-empty Zariski-open subset $U$ of $E_1/\mathfrak{m}E_1$ with the following property: for any $y \in E_1$ that has a natural image in $U$, and for any $l \geq e$ and sufficiently large $n_1, \dots, n_k > 0$ (depending on $l$), the following equalities hold:

$${\bf E}^{\bf n}M[Y]_s \cap (x^l + y^lY)M[Y]_{|{\bf n}|+s-l} = (x^l + y^lY)E_1^{n_1-l}E_2^{n_2}\cdots E_k^{n_k}M[Y]_s,$$

$$\left({\bf E}^{\bf n}M[Y]_s :_{M[Y]_{|{\bf n}|+s-l}} (x^l + y^lY)\right) \cap E_1^cE_2^{n_2}\cdots E_k^{n_k}M[Y]_{n_1-c+s} = E_1^{n_1-l}E_2^{n_2}\cdots E_k^{n_k}M[Y]_{l+s},$$
where $M[Y] = \oplus_{i \geq 0} F[Y]^i \otimes_R N.$
\end{lemma}

\begin{proof} Let $Y$ be a variable over $R$, and let $x \in E_1$. Thus, we have $x \in I_{E_1}$. According to \cite[Lemma 17.5.2]{SH}, there exist positive integers $c$ and $e$, and a non-empty Zariski-open subset $U$ of $I_{E_1}/\mathfrak{m}I_{E_1}$ with the following property: for any $y \in I_{E_1}$ that has a natural image in $U$, and for any $l \geq e$ and sufficiently large $n_1, \dots, n_k > 0$ (depending on $l$), the following equalities hold:

$${\bf I}_{{\bf E}}^{\bf n}M[Y] \cap (x^l + y^lY)M[Y] = (x^l + y^lY)I_{E_1}^{n_1-l}I_{E_2}^{n_2}\cdots I_{E_k}^{n_k}M[Y],$$

$$\left({\bf I}_{{\bf E}}^{\bf n}M[Y]:_{M[Y]} (x^l + y^lY)\right) \cap I_{E_1}^cI_{E_2}^{n_2}\cdots I_{E_k}^{n_k}M[Y] = I_{E_1}^{n_1-l}I_{E_2}^{n_2}\cdots I_{E_k}^{n_k}M[Y].$$

Moreover, by Proposition \ref{propo3.3BP}, we can choose such an open set $U$ that only consists of elements of degree one in $I_{E_1}$, meaning that it contains only elements from $E_1$. Therefore, the desired result follows by concentrating on degree $|{\bf n}|+s$ in the above equalities.
\end{proof}

\begin{proposition}\label{Proposition 17.5.1}  
Let $(R, \mathfrak{m})$ be a $d$-dimensional Noetherian local ring. Let $N$ be a finitely generated $R$-module of dimension $d$. Suppose $E_1, E_2, \ldots, E_{d+p-1}$ are $R$-submodules of a module $F$ such that $F/E_i$ has finite colength for all $i = 1, 2, \ldots, d+p-1$.  Then,  for any positive integer $l$, 
$${\rm \tilde{e}}_{\rm BR}\left(E_1, E_2, \ldots, E_{d+p-2}, E_{d+p-1}^l; N\right) = l \cdot {\rm e}_{\rm BR}\left(E_1, E_2, \ldots, E_{d+p-1}; N\right).$$
\end{proposition}
\begin{proof}
Without loss of generality, let's assume $d_k > 0$, which means that the submodules $E_i$ are repeated at least once. By using \cite[Lemma 8.4.2 and Section 8.4]{SH}, we can assume that the residue field of $R$ is infinite. According to Proposition \ref{propo3.3BP}, there exist elements $x_1, x_2, \dots, x_{d+p-1}$, where the $j$-th element is taken from $E_j$, such that they form a superficial sequence for $E_1, E_2, \dots, E_{d+p-1}$ with respect to $N$. Moreover, we can assume that for all positive integers $l$, $x_{d+p-1}^l \in E_{d+p-1}^l$ is superficial for $E_{d+p-1}^l$ with respect to $M' = S \otimes_R N / (x_1, x_2, \dots, x_{d+p-2})S \otimes_R N$.  Therefore, we have the following equalities:

\[\begin{aligned}
{\rm \tilde{e}}_{\rm BR}\left(E_1, E_2, \dots, E_{d+p-2}, E_{d+p-1}^l; N\right) &= {\rm \tilde{e}}_{\rm BR}\left(E_1^{[0]}, E_2^{[0]}, \dots, E_{d+p-2}^{[0]}, E_{d+p-1}^l; M'\right)\,\, \text{by Theorem \ref{Theorem 17.4.6degree}} \\
&= {\rm \tilde{e}}_{\rm BR}(E_{d+p-1}^l; M') \\
&= l \cdot {\rm e}_{\rm BR}(E_{d+p-1}; M') \quad \text{by Proposition \ref{Proposition 11.2.9}} \\
&= l \cdot {\rm e}_{\rm BR}(x_{d+p-1}; M') \quad \text{by Corollary \ref{Theorem 17.4.6degreecor1}} \\
&= l \cdot {\rm e}_{\rm BR}(x_1, x_2, \dots, x_{d+p-2}, x_{d+p-1}; N) \\
&= l \cdot {\rm e}_{\rm BR}(E_1, E_2, \dots, E_{d+p-1}; N) \quad \text{by Corollary \ref{Theorem 17.4.6degreecor1}}.
\end{aligned}\]

Thus, we have shown the result.
\end{proof}

\begin{lemma}\label{Lemma 17.5.4} 
Let $(R, \mathfrak{m})$ be a $d$-dimensional Noetherian local ring with an infinite residue field. Let $E_1, \dots, E_k$ be finitely generated $R$-submodules of  $F$, and $x_i \in E_i$ for $i = 1, \dots, k$. Consider a variable $Y$ over $R$. Assume that the ideals $(x_1, \dots, x_k)S$ and $I_{E_i}$ have the same height $k$ and the same radical for all $i = 1, \dots, k$. Let $\mathfrak{p}$ be a prime ideal minimal over $((x_1, \dots, x_k):_RF)$ such that ${\rm e_{BR}}((x_1, \dots, x_k)_{\mathfrak{p}};R_{\mathfrak{p}}) = {\rm e_{BR}}({E_1}_{\mathfrak{p}}, \dots, {E_k}_{\mathfrak{p}};R_{\mathfrak{p}})$. Set $B = R[Y]_{\mathfrak{p}R[Y]}$. Then, there exists a non-empty Zariski-open subset $U$ of $E_1/\mathfrak{m}E_1$ (actually of $(E_1 + (x_k))/(\mathfrak{m}E_1 + (x_k))$) that can be lifted to a non-empty Zariski-open subset $U$ of $E_1/\mathfrak{m}E_1$ such that for any preimage $y$ of an element of $U$ and for all sufficiently large integers $l$,
$${\rm e}_{\rm t}\left(\left(x_1^l + y^lY, x_2, \dots, x_k\right)_{\mathfrak{p}R[Y]}; B\right) = l \cdot {\rm e}_{BR}\left(\left(x_1, x_2, \dots, x_k\right)_{\mathfrak{p}}; R_{\mathfrak{p}}\right).$$
\end{lemma}
\begin{proof} We use induction on $k$. If $k = 1$, choose $y$ as in Lemma \ref{Lemma 17.5.2}. Then, for all sufficiently large integers $l$, $x_1^l + y^lY$ is superficial for $E_1^l(R[Y])$, hence also for $E_1^l(R[Y]_{\mathfrak{p}R[Y]})$. Thus, 
$$
\begin{array}{llll}
{\rm e_t}\left((x_1^l+y^lY)_{\mathfrak{p}R[Y]};B\right)
& = & l\cdot {\rm e_t}\left((x_1+yY)_{\mathfrak{p}R[Y]};B\right) &  \mbox{by \cite[Lemma 2.7(i)]{Bedregal-Perez3}} \\
& = & l\cdot {\rm e}_{\rm BR}\left((x_1+yY)_{\mathfrak{p}R[Y]};B\right) &  \mbox{by \cite[Theorem 4.1]{Bedregal-Perez3}} \\
& = & l\cdot {\rm {e}}_{\rm BR}\left({E_1}_{\mathfrak{p}R[Y]};B\right) &  \mbox{by Corollary \ref{Theorem 17.4.6degreecor1}} \\
& = & l\cdot {\rm e_{BR}}\left({E_1}_{\mathfrak{p}};R_{\mathfrak{p}}\right) &  \\
& = & l \cdot {\rm e_{BR}}\left((x_1)_{\mathfrak{p}};R_{\mathfrak{p}}\right) &  \mbox{by hypothesis}.
\end{array}
$$
So the case $k = 1$ is proved.

Now assume $k \geq 2$. Let $\mathfrak{p}$ be a prime ideal minimal over $\left((x_1,\dots,x_k):_RF\right)$. For $\mathfrak{q}\in {\rm Min}(R_{\mathfrak{p}}),$ set $A = R_{\mathfrak{p}}/\mathfrak{q}.$  By Lemma \ref{Lemma 17.5.3} (ii), if $\dim(A) = \dim(R_{\mathfrak{p}}),$ then $${\rm e_{BR}}\left((x_1,\cdots,x_k)_{\mathfrak{p}}(R_{\mathfrak{p}}/\mathfrak{q});A\right)\geq {\rm e_{BR}}\left({E_1}_{\mathfrak{p}}(R_{\mathfrak{p}}/\mathfrak{q}),\dots,{E_k}_{\mathfrak{p}}(R_{\mathfrak{p}}/\mathfrak{q}); A\right).$$  
Let $\Lambda=\left\{\mathfrak{q}\in \Supp(R_{\mathfrak{p}}) : \dim\left( A/\mathfrak{q}\right)=\dim(R_{\mathfrak{p}})\right\}$. By Lemmas \ref{Associativity} and \ref{AssociativityMix} follows
$$
\begin{array}{lll}  
{\rm e_{BR}}\left((x_1,\dots,x_k)_{\mathfrak{p}};R_{\mathfrak{p}}\right) 
&=&\sum\limits_{\mathfrak{q}\in \Lambda }\ell_{R_{\mathfrak{p}}}\left([R_{\mathfrak{p}}]_{\mathfrak{q}}\right) {\rm e_{BR}}\left((x_1,\dots,x_k)_{\mathfrak{p}}(R_{\mathfrak{p}}/\mathfrak{q});A\right)\\
&=&\sum\limits_{\mathfrak{q}\in \Lambda }\ell\left(R_{\mathfrak{q}}\right) {\rm e_{BR}}\left((x_1,\dots,x_k)_{\mathfrak{p}}(R_{\mathfrak{p}}/\mathfrak{q});A\right)\\
&\geq & \sum\limits_{\mathfrak{q}\in \Lambda}\ell\left(R_{\mathfrak{q}}\right) {\rm e_{BR}}({E_1}_{\mathfrak{p}}(R_{\mathfrak{p}}/\mathfrak{q}),\dots,{E_k}_{\mathfrak{p}}(R_{\mathfrak{p}}/\mathfrak{q}); A)\\
&=&{\rm e_{BR}}({E_1}_{\mathfrak{p}},\dots,{E_k}_{\mathfrak{p}}; R_{\mathfrak{p}}).
\end{array}
$$

But all the terms above must be equal (by hypothesis), so for each $A = R_{\mathfrak{p}}/\mathfrak{q}$, we have
$$
{\rm e_{BR}}\left((x_1,\dots,x_k)_{\mathfrak{p}}(R_{\mathfrak{p}}/\mathfrak{q});A\right) = {\rm e_{BR}}({E_1}_{\mathfrak{p}}(R_{\mathfrak{p}}/\mathfrak{q}),\dots,{E_k}_{\mathfrak{p}}(R_{\mathfrak{p}}/\mathfrak{q}); A).
$$

Thus, the hypotheses of the lemma hold for each $R/\mathfrak{p}$ in place of $R$, with $\mathfrak{q}$ varying over those minimal prime ideals of $\Lambda$. If the conclusion holds with $R/\mathfrak{q}$ in place of $R$, then there exists a Zariski-open non-empty subset $U_{\mathfrak{q}}$ of $E_1/\mathfrak{m}E_1$ such that the conclusion of the lemma holds for $R/\mathfrak{q}$ in place of $R$. Then, by Lemma \ref{Associativity}, the conclusion holds in $R$ for $y$ a preimage of any element of the non-empty Zariski-open subset $\cap_{\mathfrak{q}}U_{\mathfrak{q}}$ of $E_1/\mathfrak{m}E_1$. Thus, it is sufficient to prove the lemma in the case where $R_{\mathfrak{p}}$ is an integral domain.

In this case, $x_k$ is a non-zerodivisor on $S_{\mathfrak{p}}:=R_{\mathfrak{p}}[t_1,\dots,t_p]$. Set $T = S_{\mathfrak{p}}/x_kS_{\mathfrak{p}}$. Then

$$
\begin{array}{lll}
{\rm e_{BR}}\left({E_1}_{\mathfrak{p}},\dots,{E_k}_{\mathfrak{p}}; R_{\mathfrak{p}}\right) &=& {\rm e_{BR}}\left((x_1,\dots,x_k)_{\mathfrak{p}};R_{\mathfrak{p}}\right) \,\,\,\,\,\,\,\,\,  \mbox{ by assumption }\\
&=&{\rm e_{BR}}\left((x_1,\dots,x_{k-1})'_{\mathfrak{p}};T\right)\,  \mbox{by \cite[Proposition 2.6]{Bedregal-Perez2} or Remark \ref{CorteRemark}} \\
&\geq &{\rm e_{BR}}\left({E_1'}_{\mathfrak{p}},\dots,{E_{k-1}'}_{\mathfrak{p}};T\right)\,\,\,\,\,\,\,\,\,  \mbox{by Lemma \ref{Lemma 17.5.3}  }\\
&=&{\rm e_{BR}}\left((y_1,\dots,y_{k-1})'_{\mathfrak{p}};T\right) \,\,\,\,\,\,\,  \mbox{by Corollary \ref{Theorem 17.4.6degreecor1} 
for some } y_i\in E_i\\
&\geq & {\rm e_{BR}}\left({E_1}_{\mathfrak{p}},\dots,{E_k}_{\mathfrak{p}}; R_{\mathfrak{p}}\right) \,\,\,\,\,\,\, \mbox{by Lemma \ref{Lemma 17.5.3}  }
\end{array}$$
so that equality has to hold throughout. In particular,  ${\rm e_{BR}}\left((x_1,\dots,x_{k-1})'_{\mathfrak{p}};T\right)= {\rm e_{BR}}\left({E_1'}_{\mathfrak{p}},\dots,{E_{k-1}'}_{\mathfrak{p}};T\right)$. By induction on $k$, there exists a non-empty Zariski open subset $U$ of $E_1/\mathfrak{m}E_1$ such that for any preimage $y$ of an element of $U$,  for large $l$ and for large $t$, we have 
$${\rm e_{\rm t}}\left((x^l_1 + y^lY,x_2,\dots, x_{k-1})'_{\mathfrak{p}R[Y]};B'\right) = l\cdot{\rm e_{BR}}((x_1,x_2,\dots ,x_{k-1})'_{\mathfrak{p}};T),$$
where $B'=S_{\mathfrak{p}R[Y]}/x_kS_{\mathfrak{p}R[Y]}.$

Now, by \cite[Proposition 3 (ii)]{Kirby}, we finish the proof:
$$
\begin{array}{lll}
     {\rm e_{t}}\left((x^l_1 + y^lY,x_2,\dots, x_k)_{\mathfrak{p}R[Y]};B\right) &=& {\rm e_{t}}\left((x^l_1 + y^lY,x_2,\dots, x_{k-1})'_{\mathfrak{p}R[Y]};B'\right) \\
     & =&l\cdot {\rm e_{BR}}((x_1,x_2,\dots ,x_{k-1})'_{\mathfrak{p}};T)\\
     & =& l\cdot{\rm e_{BR}}\left((x_1,\dots,x_k)_{\mathfrak{p}};R_{\mathfrak{p}}\right).
\end{array}
$$
\end{proof}

\begin{lemma}\label{Lemma 17.5.5}
    Let $(R,\mathfrak{m})$ be a formally equidimensional Noetherian local ring with infinite residue field, $Y$ be a variable over $R$ and $A:=R[Y]_{\mathfrak{m}R[Y]+YR[Y]}$. Let $E_1,\dots,E_k$ be finitely generated $R$-submodules of $F$, with $x_i \in E_i$ for $i = 1,\dots, k$. Assume that the ideals $(x_1,\dots,x_k)S$ and all the $I_{E_i}$ have height k and have the same radical. Let $\Lambda$ be the set of prime ideals in $R$ minimal over $\left((x_1,\dots,x_k):_RF\right)$. Assume that for all $\mathfrak{p} \in \Lambda$, ${\rm e_{BR}}\left((x_1,\dots,x_k)_{\mathfrak{p}};R_{\mathfrak{p}}\right) = {\rm e_{BR}}\left({E_1}_{\mathfrak{p}},\dots, {E_k}_{\mathfrak{p}},R_{\mathfrak{p}}\right)$. Let $y \in E_1$ be a superficial element for $E_1,\dots,E_k$ that is not in any prime ideal minimal over $(x_2,\dots, x_k)S$. Then for all sufficiently large integers $l$, the set of prime ideals of $A$ minimal over $\left(\left[(x^l_1+y^lY,x_2,\dots,x_k)G\right]_n:_AG_n\right)$ is equal to $\{\mathfrak{p}A|\mathfrak{p}\in \Lambda\}$,  where  $G=S[Y]_{\mathfrak{m}R[Y]+R[Y]}.$ 
\end{lemma}

\begin{proof} Let $\left[(x^l_1+y^lY,x_2,\dots,x_k)G\right]_n$ be the $A$-submodule of $G_n$, for $n\geq l.$
     By the choice of $y$, the height of $I_{K_l}$ is $k$. Elements of $\Lambda$ are clearly extended to prime ideals in $A$ that are minimal over $(K_l:_AG_n)$. Suppose there exists a prime ideal $\mathfrak{q}$ in $A$, minimal over $(K_l:_AG_n)$, that is not extended from a prime ideal in $\Lambda$. By Krull's Height Theorem, \cite[Theorem B.2.1]{SH}, ${\rm ht}\mathfrak{q} \leq {\rm ht}\mathfrak{p}$. As $I_{K_l}$ has height $k$, necessarily ${\rm ht}\mathfrak{q} = {\rm ht}\mathfrak{p}$. By \cite[Lemma B.4.7]{SH}, $A$ is formally equidimensional, so that by \cite[Lemma B.4.2]{SH}, $\dim(A/\mathfrak{q}) = \dim A -{\rm ht}\mathfrak{p}$. Similarly, for each $\mathfrak{p}\in \Lambda$, $\dim(A/\mathfrak{p}) = \dim A -{\rm ht}\mathfrak{p}$.  By Additivity and Reduction Formulas, \cite[Theorem 11.2.4]{SH}, for all $t \gg 0$,
\[
\begin{array}{lll}
{\rm e}\left(\left[\frac{G}{((x^l_1+y^lY)^n,x_2^n,\dots,x_k^n)G)}\right]_t\right)&\geq & {\rm e}(A/\mathfrak{q})\ell\left(\left[\frac{G_{\mathfrak{q}}}{(x^l_1+y^lY)^n,x_2^n,\dots,x_k^n)G_{\mathfrak{q}}}\right]_t\right)+\\
&&  \sum\limits_{\mathfrak{p}\in\Lambda}{\rm e}(A/\mathfrak{p}A)\ell\left(\left[\frac{G_{\mathfrak{p}A}}{((x^l_1+y^lY)^n,x_2^n,\dots,x_k^n)G_{\mathfrak{p}A}}\right]_t\right)

\end{array}
\]

By Lechs Formula \cite[Theorem 3.1]{Bedregal-Perez3} it follows that

\[
\begin{array}{lll}
\lim\limits_{n\to\infty}\frac{{\rm e}\left(\left[\frac{G}{((x^l_1+y^lY)^n,x_2^n,\dots,x_k^n)G)}\right]_t\right)}{n^{k}}&\geq & {\rm e}(A/\mathfrak{q}){\rm e_t}({K_l}_{\mathfrak{q}};A_{\mathfrak{q}})+\sum\limits_{\mathfrak{p}\in\Lambda}{\rm e}(A/\mathfrak{p}A){\rm e_t}({K_l}_{\mathfrak{p}A};A_{\mathfrak{p}A})\\
&>& \sum\limits_{\mathfrak{p}\in\Lambda}{\rm e}(A/\mathfrak{p}A){\rm e_t}\left({K_l}_{\mathfrak{p}A};A_{\mathfrak{p}A}\right)\\
&=&  \sum\limits_{\mathfrak{p}\in\Lambda}{\rm e}(R/\mathfrak{p})\cdot l\cdot {\rm e_{BR}}\left((x_1,\dots,x_k)_{\mathfrak{p}};R_{\mathfrak{p}}\right),

\end{array}
\]

\noindent the last equality holds by Lemma \ref{Lemma 17.5.4}. By the Additivity and Reduction Formulas, \cite[Theorem 11.2.4]{SH},

\[
\begin{array}{lll}
{\rm e}\left(\left[\frac{G}{((x^l_1+y^lY)^n,x_2^n,\dots,x_k^n)G)}\right]_t\right)&\leq & {\rm e}\left(\left[\frac{G}{((x^l_1+y^lY)^n,x_2^n,\dots,x_k^n,Y)G)}\right]_t\right)\\
&=& {\rm e}\left(\left[\frac{S}{(x^{ln},x_2^n,\dots,x_k^n)S}\right]_t\right)\\
&=&\sum\limits_{\mathfrak{p}\in\Lambda}{\rm e}(R/\mathfrak{p})\ell\left(\left[\frac{S_{\mathfrak{p}}}{((x^{ln},x_2^n,\dots,x_k^n)S_{\mathfrak{p}}}\right]_t\right)

\end{array}
\]

\noindent so that by Lechs Formula \cite[Theorem 3.1]{Bedregal-Perez3} and by  \cite[Lemma 2.7(i) and Theorem 4.1]{Bedregal-Perez3}, we get
\[
\begin{array}{lll}
\lim\limits_{n\to\infty}\frac{{\rm e}\left(\left[\frac{G}{((x^l_1+y^lY)^n,x_2^n,\dots,x_k^n)G)}\right]_t\right)}{n^{k}}&\leq & \sum\limits_{\mathfrak{p}\in\Lambda}{\rm e}(R/\mathfrak{p}) {\rm e_t}\left((x^{l}_1,x_2,\dots,x_k)_{\mathfrak{p}};R_{\mathfrak{p}}\right)\\
&\leq & \sum\limits_{\mathfrak{p}\in\Lambda}{\rm e}(R/\mathfrak{p})\cdot l\cdot {\rm e_{BR}}\left((x_1,x_2,\dots,x_k)_{\mathfrak{p}};R_{\mathfrak{p}}\right)

\end{array}
\]
\noindent contradicting the earlier inequality. Thus there is no such $\mathfrak{q}$ and the result follows.
\end{proof}

\section{Main result}

In this section, we prove the main result of this paper, which is the converse of Theorem \ref{Perez-Bedregal} or the generalization of Swanson's Theorem \ref{SwansonTheorem}. In essence, our result states that if the mixed multiplicity of a family of modules is equal to the Buchsbaum-Rim multiplicity of a module generated by a sequence of elements from the family of modules, then this sequence constitutes a joint reduction of the family of modules.

\begin{theorem}[Converse of Rees' multiplicity theorem for modules] \label{17.6.1HS} Let $(R,\mathfrak{m})$ be a $d$-dimensional formally equidimensional Noetherian local ring, $E_1,\dots,E_k$ be finitely generated $R$-submodules of $F$ and $x_i\in E_i$, for $i=1,\dots,k$.  Assume that the ideals $\left(x_1,\ldots,x_k\right)S$ and $I_{E_i}$ have the same height $k$ and the same radical. If  ${\rm e_{BR}}\left((x_1,\dots,x_k)_{\mathfrak{p}};R_{\mathfrak{p}}\right) =  {\rm e_{BR}}\left({E_1}_{\mathfrak{p}},\dots, {E_k}_{\mathfrak{p}},R_{\mathfrak{p}}\right)$, for all prime ideal $\mathfrak{p}$ minimal over $\left(\left(x_1,\ldots,x_k\right):_RF\right)$, then $\left(x_1,\ldots,x_k\right)$ is a joint reduction for $\left(E_1,\dots,E_k\right).$
\end{theorem}
\begin{proof}
Let $X$ be a variable over $R$. If $(R,\mathfrak{m})$ is a local ring with finite residue field, we set $T=R[X]_{\mathfrak{m}R[X]}$ which is the localization of the polynomial ring $R$ at the prime ideal $\mathfrak{m}R[X]$. Then $T$ is a local ring with an infinite residue field and is still formally equidimensional, by \cite[Lemma B.4.7]{SH}. Furthermore, $R \subseteq T$ is a faithfully flat extension of Noetherian local rings of the same Krull dimension (see \cite[Section 8.4]{SH}), and radicals, heights, minimal prime ideals, multiplicities, and mixed multiplicities are preserved under passage to $T$; besides, some $k$-tuple of elements is a joint reduction for a $k$-tuple of submodules in $F$ if and only if it is so after passage to $T$.
Therefore, by possibly switching to $T$, we may assume that $R$ has an infinite residue field.

We will prove this theorem by induction on $k$.

When $k=1$, by hypothesis, we have $x_1S$ and $I_{E_1}$ with the same height $1$ and $\sqrt{x_1S}=\sqrt{I_{E_1}}$. Besides ${\rm e_{BR}}\left((x_1)_{\mathfrak{p}};R_{\mathfrak{p}}\right) ={\rm e_{BR}}\left({E_1}_{\mathfrak{p}},R_{\mathfrak{p}}\right),$ for all prime ideal $\mathfrak{p}$ minimal over $\left(x_1:_RF\right)$. As $R_{\mathfrak{p}}$ is formally equidimensional, by Rees's Theorem (Lemma \ref{Coro16.5.7SH}), $(x_1)_{\mathfrak{p}}$ is a reduction of ${E_1}_{\mathfrak{p}}$, so $I_{(x_1)_{\mathfrak{p}}}$ is a reduction of $I_{{E_1}_{\mathfrak{p}}}$, by Remark \ref{remarkideal}. Thus, using the definition of integral closure of $I_{{E_1}_{\mathfrak{p}}}$, we have  $I_{{E_1}_{\mathfrak{p}}}\subseteq \cap_{\mathfrak{p}}\overline{I_{(x_1)_{\mathfrak{p}}}}\cap S$, by \cite[Proposition 2.1.16]{SH}, for all prime ideal $\mathfrak{p}$ minimal over $\left(x_1:_RF\right)$. By \cite[Ratliff's Theorem  5.4.1]{SH}, $I_{E_1}\subseteq \overline{I_{(x_1)}}$, again by Remark \ref{remarkideal}, we get that $(x_1) \subseteq E_1$  is a reduction.

Suppose $k>1$. Let's reduce to the case when $R$ is a domain. Let $\Lambda$ be the set of all prime ideals in $R$ that are minimal over $\left((x_1,\dots,x_{k}):_RF\right)$. Let $\mathfrak{p}\in \Lambda$ and let $Q \in {\rm Spec}(R)$ be a minimal over $\left(\mathfrak{p}F + (x_1,\dots,x_{k}):_RF\right)$. Since $R$ is formally equidimensional, hence equidimensional and catenary, ${\rm ht}(Q)={\rm ht}(Q/\mathfrak{p}) \leq {\rm ht}({\mathfrak{p}})$, so necessarily ${\rm ht}(Q)={\rm ht}({\mathfrak{p}})$ and $Q\in \Lambda$. By Lemma \ref{Lemma 17.5.3} - item (ii),

\begin{equation}\label{inequality30}
{\rm e_{BR}}\left(\left(\left(x_1,\dots,x_{k}\right)\left(R/\mathfrak{p}\right)\right)_{Q};\left(R/\mathfrak{p}\right)_{Q}\right) \leq {\rm e_{BR}}\left(\left(E_1\left(R/\mathfrak{p}\right)\right)_{Q},\ldots,\left(E_{k}\left(R/\mathfrak{p}\right)\right)_{Q};\left(R/\mathfrak{p}\right)_{Q}\right).
\end{equation}

Then, by Associativity Formula for Buchsbaum-Rim multiplicity (Lemma \ref{Associativity}) and Reduction Formula for mixed Buchsbaum-Rim multiplicity (Lemma \ref{AssociativityMix}):
\[
\begin{array}{lll}
     {\rm e_{BR}}((x_1,\dots,x_{k})_{Q}; R_{Q})
     & = & \sum\limits_{\mathfrak{q}\in \Lambda,\, \mathfrak{q}\subseteq Q}\ell_R((R_{Q})_{\mathfrak{q}}){\rm e_{BR}}\left(\left(x_1,\dots,x_{k}\right)(R/\mathfrak{q})_{Q}; (R/\mathfrak{q})_{Q}\right) \\
     \\
     & \leq & \sum\limits_{\mathfrak{q}\in \Lambda, \mathfrak{q}\subseteq Q}\ell((R_{Q})_{\mathfrak{q}}){\rm e_{BR}}\left((E_1(R/\mathfrak{q}))_Q,\ldots,(E_{k}(R/\mathfrak{q}))_Q; (R/\mathfrak{q})_{Q}\right) \\  
     \\
     & = & {\rm e_{BR}}({E_1}_{Q},\ldots,{E_k}_{Q}; R_{Q}) \\
     \\
     & = & {\rm e_{BR}}((x_1,\dots,x_{k})_{Q}; R_{Q}),
\end{array}
\]

\noindent so, by inequality (\ref{inequality30}), follows
$${\rm e_{\rm BR}}\left(\left(\left(x_1,\dots,x_k\right)\left(R/\mathfrak{p}\right)\right)_Q; (R/\mathfrak{p})_{Q}\right) = {\rm e}_{\rm BR}\left(\left(E_1(R/\mathfrak{p})\right)_Q,\ldots,\left(E_{k})(R/\mathfrak{p})\right)_Q; (R/\mathfrak{p})_{Q}\right),$$ for each $\mathfrak{p}\in \Lambda, \mathfrak{p}\subseteq Q$. 

If the result is true for integral domains, since ${\rm ht}(Q)={\rm ht}(Q/\mathfrak{p})={\rm ht}({\mathfrak{p}})$, then $(x_1,\dots,x_k)$ is a joint reduction for $(E_1, \dots , E_k)$ with respect to $R/\mathfrak{p}$ for each $\mathfrak{p}\in \Lambda$. Then, by Proposition \ref{Proposition 1.1.5}, since the definition of joint reduction reduces to a reduction question, $(x_1,\dots,x_k)$ is a joint reduction for $(E_1, \dots, E_k)$. Thus it is sufficient to prove the theorem for integral domains.

Let $A = R[Y]_{\mathfrak{m}R[Y]+YR[Y]}$ and let $y$ be as in the statements of Lemmas \ref{Lemma 17.5.2} and \ref{Lemma 17.5.4}. Since both requirements are given by non-empty Zariski-open sets, such $y$ exists and we may choose a non-zero $y$. Thus $x^l_1 + y^lY$ is not zero for all $l$. Set $G=S[Y]_{\mathfrak{m}R[Y]+YR[Y]}$  and let $G'=G/(x^l_1 + y^lY)G$, for some large integer $l$. By Lemma \ref{Lemma 17.5.4}, if $l$ is sufficiently large,
$${\rm e_{t}}\left((x^l_1 + y^lY,x_2,\dots, x_k)_{\mathfrak{p}A};A_{\mathfrak{p}A}\right) = l\cdot {\rm e_{BR}}\left((x_1,x_2,\dots ,x_k)_{\mathfrak{p}};R_{\mathfrak{p}}\right),$$for every $\mathfrak{p}\in\Lambda$. By \cite[Proposition 3]{Kirby} and \cite[Theorem 4.1]{Bedregal-Perez3}, for $t\gg 0$,
\[ 
\begin{array}{lll}
{\rm e_{\rm t}}\left((x^l_1 + y^lY,x_2,\dots, x_k)_{\mathfrak{p}A};A_{\mathfrak{p}A}\right)&=& {\rm e_{\rm t}}\left((x_2,\dots, x_k)'_{\mathfrak{p}A};G_{\mathfrak{p}A}'\right)\\
&=& {\rm e_{BR}}\left((x_2,\dots, x_k)'_{\mathfrak{p}A};G_{\mathfrak{p}A}'\right).
\end{array}
\]

By  Lemma \ref{Lemma 17.5.2}, there exists an integer $c$ such that for all large $n_1,\dots,n_k,$

$$
\begin{array}{r}
\left({\bf E}^{\bf n}G_s:_{G_{|{\bf n}|+s-l}}(x^l+y^lY)\right)\cap E_1^{c}E_2^{n_2}\cdots E_k^{n_k}G_{n_1-c+s}=E_1^{n_1-l}E_2^{n_2}\cdots E_k^{n_k}G_{l+s}.
\end{array}$$

This, in particular, holds for all $n_1$ that are large multiples of $l$, and $c$ replaced by a larger integer that is a multiple of $l$. Thus,

\[
\begin{array}{lll}
{\rm {e}}_{\rm BR}\left({E_2'}_{\mathfrak{p}A},\ldots,{E_k'}_{\mathfrak{p}A};G_{\mathfrak{p}A}'\right)&=&{\rm \tilde{e}}_{\rm BR}\left({E_2'}_{\mathfrak{p}A},\ldots,{E_k'}_{\mathfrak{p}A};G_{\mathfrak{p}A}'\right)\\
&=&{\rm \tilde{e}}_{\rm BR}\left({E_1^l}_{\mathfrak{p}A},\ldots,{E_k}_{\mathfrak{p}A};A_{\mathfrak{p}A}\right),\,\, \mbox{by Theorem \ref{Theorem 17.4.6degree}}\\
&=&{\rm \tilde{e}}_{\rm BR}\left({E_1^l}_{\mathfrak{p}},\ldots,{E_k}_{\mathfrak{p}};R_{\mathfrak{p}}\right)\\
&=& l\cdot {\rm {e}}_{\rm BR}\left({E_1}_{\mathfrak{p}},\ldots,{E_k}_{\mathfrak{p}};R_{\mathfrak{p}}\right),\,\, \mbox{by Proposition \ref{Proposition 17.5.1}}
\end{array}
\]

By assumption and the derived equalities, we get that $${\rm {e}}_{\rm BR}\left({E_2'}_{\mathfrak{p}A},\ldots,{E_k'}_{\mathfrak{p}A};G_{\mathfrak{p}A}'\right)={\rm e_{BR}}((x_2,\dots, x_k)'_{\mathfrak{p}A};G_{\mathfrak{p}A}'), \qquad \mbox{ for all } \mathfrak{p}\in \Lambda.$$ 
\noindent  By Lemma \ref{Lemma 17.5.5}, all the minimal prime ideals over $\left([(x^l_1+y^lY,x_2,\dots,x_k)G]_n:_AG_n\right)$ are of the form $\mathfrak{p}A$, with $\mathfrak{p}\in \Lambda $.

Set $I_E = x_2I_{E_3}\cdots I_{E_k} + \cdots + x_kI_{E_2} \cdots I_{E_{k-1}}$. By \cite[Lemma B.4.7]{SH}, $G$ is local formally equidimensional, and by \cite[Proposition B.4.4]{SH}, $G'$ is formally equidimensional. By induction on $k$, $(x_2, \dots , x_k)$ is joint reduction for $(E_2,\dots , E_k)$ with respect to $G'$ wich is  equivalent to say that  $(x_2, \dots , x_k)$ is joint reduction for $(I_{E_2},\dots , I_{E_k})$ with respect to $G'$ (by Remark \ref{Remark 4.4}).

So $I_EG'$ is a reduction of $I_{E_2}\cdots I_{E_k}G'$. Thus, for sufficiently large $n$, $(I_{E_2}\cdots I_{E_k})^{n+1}G' \subseteq I_E(I_{E_2}\cdots I_{E_k})^nG'$. Hence $(I_{E_2}\cdots I_{E_k})^{n+1}G \subseteq I_E(I_{E_2}\cdots I_{E_k})^nG + (x^l_1 + y^lY )G$. By the choice of $y$ as in Lemma \ref{Lemma 17.5.2}, for possibly larger $n$, if $I_E' = I_EI_{E_1},$
$$(I_{E_1}\cdots I_{E_k})^{n+1}G\subseteq I_E' (I_{E_1}I_{E_2} \cdots I_{E_k})^nG + (x^l_1 + y^lY )I_{E_1}^{n+1-l}(I_{E_2}\cdots I_{E_k})^{n+1}G.$$ 

Thus there exists $s\in S[Y]\setminus (\mathfrak{m}R[Y]+YR[Y])$ such that
$$(I_{E_1}\cdots I_{E_k})^{n+1}S[Y] \subseteq I_E' (I_{E_1}I_{E_2} \cdots I_{E_k})^nS[Y] + (x^l_1 + y^lY )I_{E_1}^{n+1-l}(I_{E_2}\cdots I_{E_k})^{n+1}S[Y].$$

But the constant term $u$ of $s$ is a unit in $R$, so by reading off the degree zero monomials in $S[Y]$ we get
\[
\begin{array}{lll}
(I_{E_1}\cdots I_{E_k})^{n+1} &\subseteq& I_E'(I_{E_1}I_{E_2}\cdots I_{E_k})^n + x_1^l I_{E_1}^{n+1-l}(I_{E_2}\cdots I_{E_k})^{n+1}\\
&\subseteq & I_E'(I_{E_1}I_{E_2}\cdots I_{E_k})^n + x_1 I_{E_1}^{n+1}(I_{E_2}\cdots I_{E_k})^{n+1},
\end{array}
\]
which proves that $(x_1,\dots,x_k)$ is a joint reduction for $(I_{E_1},\dots, I_{E_k})$. Therefore, by Remark \ref{Remark 4.4}, we get $(x_1,\dots,x_k)$ is a joint reduction for $(E_1,\dots, E_k)$.

\end{proof}

As an immediate consequence, we have the following corollary.

\begin{corollary}
Let $(R,\mathfrak{m})$ be a $d$-dimensional formally equidimensional Noetherian local ring and $E$ be finitely generated $R$-submodules of $F$.  Assume that the ideals$\left(x_1,\ldots,x_k\right)S$ and $I_{E}$ have the same height $k$ and the same radical. If  ${\rm e_{BR}}\left((x_1,\dots,x_k)_{\mathfrak{p}};R_{\mathfrak{p}}\right) =  {\rm e_{BR}}\left({E}_{\mathfrak{p}}; R_{\mathfrak{p}}\right)$, for all prime ideal $\mathfrak{p}$ minimal over $\left(\left(x_1,\ldots,x_k\right):_RF\right)$, then $\left(x_1,\ldots,x_k\right)$ is a reduction of $R$-module $E.$    
\end{corollary}

Note that in the corollary above, we show the classic Rees' theorem for modules of finite length and the reduction criterion of B\"oger for arbitrary modules, see \cite[Corollary 16.5.7, Theorem 16.5.8]{SH} respectively.

\end{document}